\documentclass[12pt]{article}
\usepackage{amsmath}
\usepackage{amsthm}
\usepackage{amssymb}
\usepackage{float}

\newtheorem{theorem}{Theorem}
\newtheorem{lemma}[theorem]{Lemma}
\newtheorem{proposition}[theorem]{Proposition}
\newtheorem{Corollary}[theorem]{Corollary}
\newtheorem{hypo}{Hypothesis}
\newtheorem{notation}{Notation}

\newtheorem{rem}{Remark\!\!}

\newcommand{\noi}{\noindent}
\newcommand{\non}{\nonumber}
\newcommand{\disp}{\displaystyle}

\def\eps{\varepsilon}
\def\un{{\bf 1}}

\def\d{\partial}
\def\N{{\mathbb N}}
\def\R{{\mathbb R}}
\def\C{{\mathbb C}}
\def\Z{{\mathbb Z}}

\def\rhoN{\rho^N}
\def\rhond{\rho_{\rm od}}
\def\rhod{\rho_{\rm d}}
\DeclareMathOperator{\im}{{\rm Im}}
\DeclareMathOperator{\re}{{\rm Re}}
\DeclareMathOperator{\range}{{\rm Range}}
\DeclareMathOperator{\Ker}{{\rm Ker}}
\def\Psiave{\langle\Psi_\eps\rangle}
\def\Psiavezero{\langle\Psi_0\rangle}
\def\Psiosc{\Psi_\eps^{\rm osc}}
\def\Psidom{\langle\Psi_\eps\rangle^{\rm dom}}
\def\PsidomN{\langle\Psi_\eps\rangle^{{\rm dom},N}}
\def\Psineg{\langle\Psi_\eps\rangle^{\rm neg}}
\def\Psiapp{\Psi^{\rm app}}
\def\Psising{\langle\Psi_\eps\rangle^{\rm sing}}
\def\Psinonsing{\langle\Psi_\eps\rangle^{\rm nonsing}}

\def\Wmod{{W_\eps^{\rm mod}}}
\def\Wout{W^{\rm pol}}

\def\rhoun{\rhod^{(1)}}
\def\rhodeux{\rhod^{(2)}}
\def\rhoapp{\rhod^{\rm app}}

\def\rhoout{\rhod^{\rm pol}}
\def\rhoNdeux{\rhod^{N,(2)}}
\def\rhoe{\underline{\rho}}

\def\op{_\sharp}
\def\transp{\!~^{t}\!}
\def\suchthat{;\ }

\newcounter{saveeqn}
\newcommand{\alpheqn}{%
  \setcounter{saveeqn}{\value{equation}}%
  \stepcounter{saveeqn}%
  \setcounter{equation}{0}%
  \renewcommand{\theequation}{\arabic{saveeqn}.\alph{equation}}%
}
\newcommand{\reseteqn}{%
  \setcounter{equation}{\value{saveeqn}}%
  \renewcommand{\theequation}{\arabic{equation}}%
}

\title{From Bloch model to the rate equations II: the case of almost
degenerate energy levels}
\author{B. Bid\'egaray-Fesquet $^{(1)}$, F. Castella $^{(2)}$,
E. Dumas $^{(3)}$ and M. Gisclon $^{(4)}$}
\date{}
\begin{document}

\maketitle

\begin{center}
\small
(1) LMC - IMAG, UMR 5523 (CNRS-UJF-INPG) \\
B.P. 53, 38041 Grenoble Cedex 9 - France\\
email: Brigitte.Bidegaray@imag.fr
\end{center}

\begin{center}
\small
(2) IRMAR, UMR 6625 (CNRS-UR1) \\
Universit\'e de Rennes 1\\
Campus de Beaulieu, 35042 Rennes Cedex - France\\
email: Francois.Castella@univ-rennes1.fr
\end{center}

\begin{center}
\small
(3) Institut Fourier, UMR 5582 (CNRS-UJF) \\
100 rue des Math\'ematiques \\
Domaine Universitaire \\
BP 74, 38402 Saint Martin d'H\`eres - France \\
email: edumas@ujf-grenoble.fr
\end{center}
\begin{center}
\small
(4) LAMA, UMR 5127 (CNRS - Universit\'e de Savoie) \\
UFR SFA, Campus Scientifique, \\
73376 Le Bourget-du-Lac Cedex - France \\
email: gisclon@univ-savoie.fr
\end{center}

\begin{abstract}
Bloch equations give a quantum description of the coupling between an atom and 
a driving electric force. In this article, we address the asymptotics of these 
equations for high frequency electric fields, in a weakly coupled regime. We 
prove the convergence towards rate equations (\textit{i.e.} linear Boltzmann
equations, describing the transitions between energy levels of the atom). We 
give an explicit form for the transition rates. 

This has already been performed in \cite{BCD} in the case when the energy 
levels are fixed, and for different classes of electric fields: quasi or 
almost periodic, KBM, or with continuous spectrum. Here, we extend the study 
to the case when energy levels are possibly almost degenerate. However, we 
need to restrict to quasiperiodic forcings. The techniques used stem from 
manipulations on the density matrix and the averaging theory for ordinary 
differential equations. Possibly perturbed small divisor estimates play a key
r\^ole in the analysis. 

In the case of a finite number of energy levels, we also precisely analyze the 
initial time-layer in the rate aquation, as well as the long-time convergence 
towards equilibrium. We give hints and counterexamples in the infinite 
dimensional case.  
\end{abstract}


\section{Introduction}

Bloch equations model the time evolution of a quantum mechanical system
described in the density matrix formulation and driven by an electromagnetic
field. This formalism is very precise but sometimes difficult to interpret
and to use in practical simulations. It is therefore useful to find asymptotic
models under appropriate scaling assumptions. In a former article \cite{BCD} 
such a program is performed and leads to the rigorous derivation of rate 
equations, which are often used in the physics literature. Here, we want to 
extend these results to the case when the energy levels of the system are 
almost degenerate. There are many examples of such almost degeneracies. This 
is the case for example of Zeemann hyperfine structures in complex molecules, 
or quantum dots submitted to an external magnetic field. High levels of an atom
are also almost degenerate since there is an infinite number of levels with
accumulation value at the ionisation energy.

\subsection{Bloch equations}

According to the quantum theory, matter is described \emph{via} a density 
matrix $\rho$, whose diagonal entry $\rho(t,n,n)$ is --in the eigenstates 
basis-- the population of the $n$-th energy level at time $t$, and the 
off-diagonal entry $\rho(t,n,m)$ is linked to the transition probability from 
level $n$ to level $m$ (conditioned by the corresponding populations). One may 
think of a collection of identical, uncoupled atoms, with discrete energy 
levels. We refer the reader to \cite{Boh,Boy,CTDRG,Lo,NM,SSL,Bi} for textbooks 
about wave/matter interaction issues, where Bloch equations occur. To treat 
the mathematical problem, we use a dimensionless version of these equations 
and consider that the density matrix $\rho(t,n,m)$ is governed by:
\begin{eqnarray}
\label{bloch}
\eps^2 \d_t \rho(t,n,m)
& = & -i \omega_\eps(n,m) \rho(t,n,m) + Q_\eps(\rho)(n,m)
\\
\non
&& + i \eps \sum_k
\left[ {\cal V}\left(\frac{t}{\eps^2},n,k\right) \rho(t,k,m)
- {\cal V}\left(\frac{t}{\eps^2},k,m\right) \rho(t,n,k) \right].
\end{eqnarray}
Integers $n \in \N$, $m \in \N$, and $k \in \N$ are labelling discrete energy
levels, and time $t$ belongs to $\R^+$. In the case of a finite number of
energy levels, we add the restriction $n \leq N$, $m \leq N$, $k \leq N$.

The dependence of the density matrix on the small parameter $\eps$
will always be implicit. Here, the system is
forced by a high frequency electromagnetic wave which is described by
its oscillatory amplitude $\phi(t/\eps^2)$ and contributes to the quantity
\[
{\cal V}\left(\frac{t}{\eps^2},n,m\right) 
= \phi\left(\frac{t}{\eps^2}\right) V(n,m),
\]
where the interaction coefficient $V(n,m)\in\C$ is (up to a rescaling) an
entry in the dipole moment matrix which is Hermitian: $V(m,n)=V(n,m)^*$.

The parameter $\eps$ occurring in the scaling plays two r\^oles. First, the 
coupling between atoms and the wave is small, of order $\eps$, and its 
cumulated effects are considered over long time scales, of size $1/\eps^2$. 
This is the setting of the so-called weak coupling regime (see 
\cite{Sp2,Sp3,VH1,VH2}). Second, the scaling in Eq. \eqref{bloch} produces 
resonances between the eigenfrequencies $\omega_k/\eps$ of the free atom and 
the electromagnetic source $\phi(t/\eps^2)$.

The quantity $\omega_\eps(n,m) = \omega(n,m) + \eps^p \delta(n,m)$ is the 
transition energy between levels $n$ and $m$. In the non-degenerate framework
$\omega(n,m) = \omega(n) - \omega(m) \in\R$ is the difference between the
energies $\omega(n)$ and $\omega(m)$ of levels $n$ and $m$ respectively. Here 
we model almost degenerate levels replacing the energy $\omega(n)$ by
$\omega(n)+\eps^p\delta(n)$ where $p>0$ and $\delta(n)\in\R$. Therefore
$\delta(n,m)=\delta(n)-\delta(m)$. We notice that a large value for $p$ means
that levels $n$ and $m$ are very close from one another if $\omega(n,m)=0$.

Last, relaxation terms are modelled by the operator $Q_\eps$ which reads
\begin{align*}
Q_\eps(\rho)(n,m) & = - \eps^\mu\gamma(n,m) \rho(t,n,m), 
& \textrm{ if }n\neq m, \\
& = \eps^2 \left( 
\sum_k \left[ W(k,n) \rho(t,k,k) - W(n,k) \rho(t,n,n) \right]
\right), & \textrm{ if }n = m.
\end{align*}
We assume that the Pauli coefficient $W(n,m)$ is non-negative, as well
as the longitudinal coefficient $\gamma(n,m)$, which is besides symmetric:
$\gamma(n,m)=\gamma(m,n)$. To simplify further notations, we extend its
definition to the case when $n=m$ introducing $\gamma(n,n)=0$. 

We assume that $0\leq\mu<1/2$. This is already the case in the
nondegenerate case
\cite{BCD}. This threshold value for $\mu$ arises in the estimates and we do
not know whether it is optimal or not. However it is no wonder that there is a
threshold value. Indeed, in the case when $\mu=0$, the initial Bloch equation
\eqref{bloch} is time-irreversible and the asymptotic equation that we derive
in this paper is also time-irreversible. On the other hand, in the  opposite
case when every coefficient $\gamma(n,m)$ is identically zero (which can also 
be interpreted as $\mu=\infty$), the initial Bloch equation is time-reversible
and the nature of the problem has changed.

The relaxation operator $Q_\eps$ determines an equilibrium state, to which the
system has a tendency to relax when no electromagnetic wave is applied. If the
threshold $\mu<1/2$ is not necessarily optimal, it is very important that
$\mu<2$. This means that the off-diagonal terms of the density matrix are
more rapidly decreasing that the diagonal terms. We therefore expect them to
play no significant r\^ole on the time scale $1/\eps^2$.

Finally, we associate to Eq.~\eqref{bloch} an initial datum $\rho(0,n,m)$
which satisfies
\begin{equation}
\label{datum}
\rho(0,n,m)=0, \quad \forall n \neq m,
\qquad \textrm{ and } \qquad
\rho(0,n,n)\geq0, \quad \forall n.
\end{equation}
The vanishing of the off-diagonal entries means that we are dealing with 
``well prepared'' initial data. This is a standard assumption in the field 
(see \textit{e.g.} \cite{KL1,KL2,Zw}).

\begin{rem}
Note that the small coupling in Eq. \eqref{bloch} \textrm{a priori} induces 
perturbations of size $1/\eps$ on time scales $1/\eps^2$. As we shall see, the 
very Hamiltonian nature of the equation actually makes these perturbations of 
size $1$.
\end{rem}

\subsection{Rate equations}

To describe the asymptotic dynamics, we are only interested in the diagonal
part of the density matrix (the populations) and set
\[
\rhod(t,n):=\rho(t,n,n),
\]
which is the occupation number of the $n$-th level. The limiting model we look
for as $\eps\to0$ is a system of rate equations which are Boltzmann type
equations, in the form
\begin{equation}
\label{rate}
\d_t \rhoapp = (\Wmod)\op \rhoapp,
\end{equation}
where the approximate populations $\rhoapp$ are viewed as a vector, and the 
modified relaxation matrix operator $(\Wmod)\op$ is defined from coefficients 
$\Wmod(n,k)$ \textit{via} the following notation.

\begin{notation}
\label{Not_op}
To coefficients $A(n,k)$, $n \neq k$ (which may possibly be time dependent: 
$A(t,n,k)$) we associate the matrix $A\op$ (respectively $A\op(t)$) through
\begin{align*}
A\op(n,k) & := A(k,n) && \textrm{if } n \neq k, \\
A\op(n,n) & := -\sum_{m \neq n} A(n,m) && \textrm{else}.
\end{align*}
If $A(n,k) \in l_n^\infty l_k^1\cap l_k^\infty l_n^1$, $A\op$ is a bounded 
operator on $l^p$, $1 \leq p \leq \infty$ (see Appendix~\ref{Sec_relax}).
\end{notation}

Apart from \cite{BCD}, in the past few years an extensive attention has been
paid on the rigourous derivation of Boltzmann type equations from dynamical
models of (classical or quantum) particles or models for the interaction of
waves with random media. Convergence results in the case of an electron in a
periodic box may be found in \cite{Ca1,Ca2,Ca3}. We also mention the
non-convergence result established in \cite{CP1,CP2} in a particular, periodic
situation. For the case when an electron is weakly coupled to random
obstacles, the reader may refer to \cite{EY,Sp1,Sp2,Sp3} and \cite{KPR} for
the formal analysis performed. The computation of the relevant cross-sections
is performed in \cite{Ni}. All these results address the case of a linear
Boltzmann equation. A nonlinear case is studied in \cite{BCEP}.

As proved in \cite{BCD}, in the non-perturbed case ($\delta(n)\equiv0$) the
action of the wave transforms the relaxation operator $W\op$ associated with 
matrix $W$ into a limit relaxation operator $(\Wmod)\op$.

The goal of this article is to specify the precise form of matrix $\Wmod(n,k)$
for the perturbed case for a fixed $\eps$. Once this rate equation is
derived, we also give results about the limit $\eps\to0$ and the size of the
time-layer leading to equilibrium according to the value of $\mu/p$.

\subsection{Main results}

\subsubsection{The modified relation operator $(\Wmod)\op$}

The modified transition rates $\Wmod(n,k)$ is defined by
\[
\Wmod(n,k) := \Psidom(n,k) + W(n,k),
\]
and the non-negative coefficient $\Psidom(n,k)$ is given (see
Proposition~\ref{Pr_relax} and Theorem~\ref{Th_relax}) by 
\begin{equation*}
\begin{aligned}
\Psidom(n,k) & := 2 |V(n,k)|^2 \sum_{ \beta \in \Z^r}
\frac{\gamma(k,n)} {\gamma(k,n)^2 + |\omega(n,k)+
\beta\cdot\omega+\eps^p \delta(k,n)|^2}
|\phi_\beta|^2 && \textrm{ if } \mu=0, \\
& := 
2 |V(n,k)|^2 \frac{\eps^\mu \gamma(k,n)}
{\eps^{2\mu}\gamma(k,n)^2 + \eps^{2 p} \delta(k,n)^2}
\sum_{ \beta \in \Z^r \suchthat \omega(k,n)+\beta\cdot\omega=0} |\phi_\beta|^2 
&& \textrm{ if } \mu>0.
\end{aligned}
\end{equation*}
Here, the electromagnetic wave $\phi$ is determined by its frequency vector
$\omega\in\R^r$ and its Fourier coefficients $\phi_\alpha$:
\[
\phi(t)=\sum_{\alpha \in \Z^r} \phi_\alpha \exp(i \alpha \cdot \omega t),
\textrm{ where we denote } 
\alpha \cdot \omega:=\alpha_1 \omega_1 + \cdots + \alpha_r \omega_r .
\]

Under assumptions that we specify in Section~\ref{Sec_fonctionnel}),
we prove the following: if $\rho$ is a solution to 
Eq.~\eqref{bloch}, and $\rhoapp$ the solution to Eq.~\eqref{rate} with the 
same initial datum, then for all $T>0$, there exists $C>0$ such that
\[
\|\rhond\|_{L^\infty([0,T],l^1)} \leq C \eps^{1-\mu}\qquad
\textrm{ and } 
\|\rhod-\rhoapp\|_{L^\infty([0,T],l^2)} \leq C (\eps^\mu + \eps^{1-2\mu}).
\]

\subsubsection{The asymptotic state}

We study the convergence of the solution $\rho$ to Bloch 
equation~\eqref{bloch} on some \textit{fixed} time interval $[0,T]$, as 
$\eps$ goes to zero. This dynamics is thus given by the corresponding 
solution $\rhoapp$ to Eq.~\eqref{rate}. Since the operator $(\Wmod)\op$ has a 
non-positive spectrum (see the Appendix~\ref{Sec_relax}), one may think that, 
as time grows, $\rhoapp$ approaches some equilibrium state, \textit{i.e.} a 
stationary state $\rhoe \in l^2$ belonging to the kernel of $(\Wmod)\op$. We 
describe this convergence carefully in the case of a finite number of quantum 
levels. In the case of an infinite number of levels, we exhibit simple 
examples for which convergence does not occur.

\subsubsection{The case of a finite number of levels}
\label{finitN}

The evolution of $\rhod$ can be summarized as:
\begin{itemize}
\item[(i)] The levels $n$ which do not resonate with the wave 
($\Psidom(n,k)=0$ for all $k$) and are not coupled to others \textit{via} 
relaxations either ($W(n,k)=0$ for all $k$) remain constant.
\item[(ii)] Polarization: the remaining levels belonging to the space 
($l^2$)-orthogonal to the kernel of the singular part 
$\Psising\op$ of $\Psidom\op$ vanish in time $O(\eps^\sigma)$. The
precise definitions of $\Psising\op$ and $\sigma>0$ are given in
Section~\ref{Sec_lim}.  
\item[(iii)] After this time, the evolution of the remaining levels is given 
by an $\eps$-independent system for $\Pi\rhoapp$,
\[
\partial_t \Pi \rhoapp = \Pi (W + \Psi_0^{\rm nonsing})\op \Pi \rhoapp.
\]
Here, $\Pi$ is the orthogonal projector onto the space orthogonal to the ones 
treated in the first two steps. Also, we set 
$\Psinonsing := \Psidom - \Psising$, and $\Psi_0^{\rm nonsing}$ is the value 
of $\Psinonsing$ as $\eps$ vanishes. For each initial datum $\rho(0)$, the 
solution to this system converges, as $t$ goes to infinity, to some (unique) 
equilibrium state $\rhoe$ (see Appendix~\ref{Sec_relax}).
\end{itemize}

In the unperturbed case $\delta=0$ \cite{BCD}, we have $\sigma=\mu$. The 
perturbation can affect the exponent $\sigma$, the transition rates given by 
$\Psising\op$, as well as the final equilibrium (determined by 
$\Psinonsing\op$). The precise value of $\sigma$ and the possibility of 
modifying the limiting operator according to the value of $\mu/p$ are given in 
Table~\ref{table1}. 

\begin{table}[H]
\[
\begin{array}{c|c|p{7cm}}
\mu/p & \sigma & difference with the unperturbed case \cite{BCD} \\
\hline \hline 0 \leq \mu/p < 1 & \mu & none \\
\hline \mu/p = 1 & \mu & transition rates \\
\hline 1 < \mu/p < 2 & 2p-\mu & transition rates and time-layer \\
\hline 2 \leq \mu/p < \infty & \mu & projector $\Pi$ and asymptotic state
\end{array}
\]
\caption{\label{table1}Consequences of the perturbation for a finite 
number of levels.}
\end{table}

\subsubsection{The case of an infinite number of levels}

For an infinite number of levels, Proposition~\ref{Pr_relax} and
Theorem~\ref{Th_relax} hold. They give the approximation of the solution to 
Eq.~\eqref{bloch} by the solution to Eq.~\eqref{rate}. The transition
rates $W(n,k)$ are again transformed into $\Wmod(n,k)$. However, the analysis 
of the asymptotic behavior of the solution to Eq.~\eqref{rate} is more 
intricate than in the finite case. In Appendix~\ref{Sec_astate}, we give 
examples for which there \textit{does not exist} any possible equilibrium 
state in $l^2$ (see Corollary~\ref{noeq}).

It is still true that non-interacting levels remain constant, as in item (i) 
of Section~\ref{finitN}. Nevertheless, we are unable to prove the
precise polarization property of the two other items, even if we
expect that the results of Table~\ref{table1} hold. An equivalent to 
Table~\ref{table1} for an infinite number of levels is given in 
Section~\ref{Sec_5infini}, for the values of $\mu/p$ for which we can conclude.

\subsection{Outline}
\label{Sec_outline}

The article is organized as follows. Section~\ref{Sec_fonctionnel} is devoted 
to the introduction of the precise notations and assumptions needed in the 
sequel. In Section~\ref{Sec_pop} a \textit{closed} equation which governs 
populations only is derived from the Bloch equations up to an approximation 
of order $O(\eps^{2(1-\mu)})$. This equation is already a Boltzmann type
equation, but the transition rates depend on time and $\eps$. The 
off-diagonal terms are proven to be negligible, of order $O(\eps^{1-\mu})$. 
This transformation uses classical arguments for the Bloch equation in the 
weak coupling regime (see \cite{Ca1,Ca2,Ca3} and also \cite{KL1,KL2,Kr,Zw} 
for this point).

In Section~\ref{Sec_ave} a new Boltzmann type equation is derived removing 
the time-dependence of the transition rates up to an approximation of order 
$O(\eps^{1-2\mu})$. These transition rates still include some terms 
which do not correspond to a resonance between the wave and the matter. These 
contributions are removed in Section~\ref{Sec_dom} with a new error of order
$O(\eps^\mu)$. The thus obtained equation is the rate equation we want to 
derive in this article. The main theorem in this Section (namely 
Theorem~\ref{Th_relax}) yields the form of the transition rates in the rate 
equation. The techniques used are those of the averaging theory for ordinary 
differential equations (see \cite{LM,SV}). Diophantine estimates play 
naturally a key r\^ole in the analysis. Lemma~\ref{diophstable} ensures that 
although the Diophantine condition is not stable with respect to small
perturbations, violations of the condition only occur for large values of the 
indices, that are compensated by extra smoothness assumptions (see
\cite{CCCFLLT} for a similar argument in another context).

The sequel of the article addresses the analysis of the limit process 
$\eps\to0$ in this equation. Section~\ref{Sec_lim} is devoted to the case 
study according to the value of $\mu/p$ leading to Table~\ref{table1} in the 
finite dimensional case (finite number of levels) and to partial results in 
the infinite dimensional case (see Table~\ref{table2}). The justification of 
the restriction to a finite number of levels is given in 
Section~\ref{Sec_fini}: we show  that there exists a number $N$ such that 
considering only the $N$ first levels implies an additional error of order 
$o(1)$. The number $N$ \textit{a priori} depends on $\eps$, except when 
resonances do not occur beyond some level.

Finally, we give in Section~\ref{Sec_app} the proofs of several lemmas, 
concerning continuity and non-positiveness of the relaxation operators first, 
then about existence and uniqueness of the associated equilibrium state, and 
finally, we show the genericity of the Diophantine condition~\ref{dioph} (see 
below).

\section{Functional setting}
\label{Sec_fonctionnel}

The choice for functional spaces is first guided by some physical properties
of the density matrix. Throughout this article the initial datum is taken such
that
\begin{equation}
\label{initial}
\rho(0,n,m)=0, \textrm{ if } n\neq m, \qquad \rho(0,n,n)\geq0 \qquad
\textrm{ and } \qquad \sum_n \rho(0,n,n) < \infty.
\end{equation}
The energy levels are assumed to be bounded:
\[
\left(\omega(n)\right)_{n \in \N} \in l^\infty
\textrm{ and }
\left(\delta(n)\right)_{n \in \N} \in l^\infty
\]
(which is physically relevant, since these energies are bounded by the
ionisation energy), as well as the relaxation coefficients
\[
\left(\gamma(n,m)\right)_{n,m \in \N} \in l^\infty, \qquad
\inf_{n\neq m} \gamma(n,m):=\gamma >0.
\]
We also suppose that the Pauli coefficients saisfy
\[
\|W\|_{l^\infty_k l^1_n \cap l^\infty_n l^1_k}
:= \sup_k \sum_n |W(n,k)| + \sup_n \sum_k |W(n,k)| < \infty.
\]
Matrix $W$ describes the relaxation to a thermodynamic equilibrium. Entries 
$W(n,m)$ and $W(m,n)$ are therefore related by the standard
microreversibility relation (see \textit{e.g.} \cite{BBR,Bi})
\begin{equation} 
\label{thermo} 
W(n,m) = \exp\left(\frac{\omega(m,n)}{T}\right)W(m,n),
\end{equation} 
where $T$ is a normalized temperature. This specific form is of great 
importance when describing the equilibrium states in Sections~\ref{Sec_lim} and
\ref{Sec_relax}. 

We recall that the interaction matrix $V$ is Hermitian. Thus we simply assume
that
\[
\|V\|_{l^\infty l^1} := \sup_k \sum_n |V(n,k)| < \infty.
\]

Classical ordinary differential equation arguments (see \textit{e.g.}
\cite{Ca1}) allow to state the existence and uniqueness of solutions to 
System~\eqref{bloch} for initial data in $l^1$. Indeed, since we assumed that 
${\cal V}(t,n,m)$ belongs to 
$L^\infty(\R^+, l^\infty_m l^1_n \cap l^\infty_n l^1_m)$, and also that 
$W(n,m) \in l^\infty_m l^1_n \cap l^\infty_n l^1_m$, we may use the estimate 
of Lemma~\ref{schur}. The operators involved in Eq.~\eqref{bloch} are thus 
continuous on $L^\infty(\R^+,l^1)$. This implies in a straightforward way 
that these solutions exist and have the following regularity: 
$\rho \in C^0(\R^+,l^1)$ and $\d_t \rho \in L^\infty(\R^+,l^1)$, for each 
$\eps>0$.

A weaker result is the existence and uniqueness in $l^2$ and we will have to 
restrict to this case from Section~\ref{Sec_ave} on.

The summation and positiveness properties are preserved through the time 
evolution if the density matrix is solution to the Bloch equations 
\eqref{bloch}. More precisely we have the following lemma (points $(ii)$, 
$(iii)$, $(iv)$ are addressed in \cite{Li,BBR,Ca3}).

\begin{lemma}
\label{conservation}
Let $\rho(t=0)$ satisfy conditions \eqref{initial}. Then, under the above 
assumptions, for all $t \in \R$,
\begin{itemize}
\item[(i)] there exists a unique solution $\rho\in{\cal C}^0(\R^+,l^1)$ to
Eq.~\eqref{bloch}.
\item[(ii)] $\rho(t)$ is Hermitian: $\rho(t,n,m)=\rho(t,m,n)^*$,
\item[(iii)] the trace of $\rho$ is conserved:
\begin{equation}
\label{lindblad}
\sum_n \rho(t,n,n)=\sum_n \rho(0,n,n)<\infty,
\end{equation}
\item[(iv)] positiveness of populations is conserved: $\rho(t,n,n) \geq 0$.
\end{itemize}
\end{lemma}

We stress the importance of items \textit{(iii)} and \textit{(iv)}, first 
established in \cite{Li}. They give a nontrivial $l^1$ estimate for the 
diagonal part $\rhod$. This proves to be crucial in 
Proposition~\ref{Pr_rhodun} (see also \cite{Ca3} for a situation where the 
oscillations are much more difficult to handle). 

The electromagnetic wave is real and bounded: $\phi(t) \in \R$ and
$\phi \in L^\infty(\R)$. Although a larger class of high frequency waves is
addressed in \cite{BCD}, we have to restrict here to the quasiperiodic
case. The amplitude $\phi$ is characterized \textit{via} its frequency vector
$\omega\in\R^r$ and its Fourier coefficients $\phi_\alpha$:
\begin{equation}
\label{Fourier}
\phi(t)=\sum_{\alpha \in \Z^r} \phi_\alpha \exp(i \alpha \cdot \omega t)
\; \textrm{ where } \;
\phi_\alpha^*=\phi_{-\alpha} \;
\; \textrm{ and } \;
\sum_{\alpha \in \Z^r} |\phi_\alpha|^2 < \infty.
\end{equation}
Here we denote 
$\alpha \cdot \omega:=\alpha_1 \omega_1 + \cdots + \alpha_r \omega_r$. \\

We need a certain number of assumptions on the wave and the interaction
coefficients. Those already stressed will be valid (and not recalled) in the
whole article. However, in some places, we will have to reinforce the decay
assumptions and therefore we also label some specific hypotheses.

As it is usual in the field of oscillations in ordinary differential equations
and averaging techniques (see \cite{Ar,SV,LM}), we introduce a Diophantine
condition on the frequency vector $\omega$. 

\begin{hypo}[Diophantine condition]
\label{dioph}
There exists a (small) number $\eta>0$, and a constant $C_\eta>0$, such that
\alpheqn
\begin{eqnarray}
&&
\non
\forall \alpha=(\alpha_1,\ldots,\alpha_r) \in \Z^r \setminus \{ 0 \}, \quad
\forall (n,k) \in \N^2 \text{ such that }
\alpha \cdot \omega + \omega(n,k) \neq 0,
\\
&&
\label{diopha}
\qquad
|\alpha \cdot \omega+\omega(n,k)| \geq
\frac{C_\eta}{(1+|\alpha|)^{r-1+\eta} (1+n)^{1+\eta} (1+k)^{1+\eta}},
\end{eqnarray}
and similarly
\begin{equation}
\label{diophb}
\forall
\alpha \in \Z^r\setminus\{0\},
\qquad
|\alpha \cdot \omega| \geq
\frac{C_\eta}{(1+|\alpha|)^{r-1+\eta}}.
\end{equation}
\reseteqn
\end{hypo}

\begin{rem}
Given once and for all a fixed $\eta>0$, we can classically claim (see 
\cite{Ar}) that there exists a constant $C_\eta>0$, depending on 
$\omega(n,m)$ and on $\eta$, such that for almost all value of the frequency 
vector $\omega=(\omega_1,\ldots,\omega_r)$ Hypothesis~\ref{dioph} is 
satisfied. This is proved in Appendix~\ref{Sec_dioph}. This condition is 
therefore not much restrictive.
\end{rem}

The same kind of estimate is also needed in the case when $\alpha=0$, which
means that the energies do not converge too fast towards the ionisation energy.

\begin{hypo}[Convergence towards the ionisation energy] 
\label{speed}
There exists a constant $C_\eta>0$, such that
\[
\forall (n,k) \in \N^2, \quad |\omega(n,k)| \geq
\frac{C_\eta}{(1+n)^{1+\eta} (1+k)^{1+\eta}},
\]
where $\eta$ is the number occurring in the Diophantine assumption
(Hypothesis~\ref{dioph}).
\end{hypo} 

From Section~\ref{Sec_ave} on, we impose the following hypothesis in
order to average transition rates in time.

\begin{hypo}[Smoothness assumption]
\label{smooth}
The Fourier coefficients $\phi_\alpha$ satisfy
\[
\sum_\alpha (1+|\alpha|)^{r-1+\eta} |\phi_\alpha|^2 < \infty,
\]
where $\eta$ is the number occurring in the Diophantine assumption
(Hypothesis~\ref{dioph}).
\end{hypo}

From Section~\ref{Sec_dom} on, we need the following two hypotheses to justify
the sorting out of resonant and non resonant contributions.

\begin{hypo}[Reinforced smoothness assumption]
\label{extrasmooth}
There exists $N_\eta>2\mu/p$ such that the Fourier coefficients $\phi_\alpha$
satisfy
\[
\sum_\alpha (1+|\alpha|)^{(r-1+\eta)N_\eta} |\phi_\alpha|^2 < \infty,
\]
where $\eta$ is the number occurring in the Diophantine assumption
(Hypothesis~\ref{dioph}).
\end{hypo}

\begin{hypo}["Far from continuum" assumption]
\label{far}
The interaction coefficients satisfy
\[
\sup_n \sum_m \left( (1+n)^{1+\eta} (1+m)^{1+\eta} \right)^{N_\eta} |V(n,m)|^2
< \infty,
\]
where $\eta$ is the number occurring in the Diophantine assumption
(Hypothesis~\ref{dioph}) and $N_\eta$ is given by 
Hypothesis~\ref{extrasmooth}.
\end{hypo}

This last hypothesis means that only low levels (\textit{i.e.} levels which 
are far from the continuous spectrum or ionisation threshold) really interact 
with the wave with a significant contribution. To restrict to a finite number 
of levels we also use an assumption on the interaction of low and high levels
\textit{via} the relaxation operator.

\begin{hypo}[Weak interaction of low and high energy levels]
\label{Wdecay}
The longitudinal relaxation coefficients satisfy
\[
\sup_n \sum_m (1+n)^K (1+m)^K |W(n,m)| < \infty,
\]
for some $K>0$.
\end{hypo}

We introduce a last notation useful for describing approximations
of $\rho$.

\begin{notation}
For $s\geq0$, and $q\geq1$, the symbol 
$O_{L^\infty_{\rm loc}(\R,l^q)}\left(\eps^s\right)$ means that for all $T>0$, 
there exists a constant $C>0$, that does not depend on $\eps$, such that the
corresponding term is bounded:
\[
\|O_{L^\infty([0,T],l^q)}(\eps^s)\|_{L^\infty([0,T],l^q)} \leq C \eps^s.
\]
\end{notation}

In the whole article $C$ denotes constants which do not depend on the (small) 
parameter $\eps$. It however possibly depends on all the coefficients of the 
problem and on the initial data, but we will never make this dependence 
explicit.

\section{A first closed equation for populations}
\label{Sec_pop}

\subsection{Populations and coherences}

In the same way as we defined the populations $\rhod(t,n)$, we define the
coherences as the off-diagonal part of the density matrix
\[
\rhond(t,n,m)=\rho(t,n,m) \un[n \neq m].
\]
Set
\[
\Omega^\eps(n,m) 
:= -i \omega(n,m)- i \eps^p \delta(n,m) - \eps^\mu \gamma(n,m),
\]
where we stress that $\Omega^\eps(n,n)=0$. With these notations, 
Eq.~(\ref{bloch}) reads for the populations:
\alpheqn
\begin{eqnarray}
\label{eqd}
\d_t \rhod(t,n) & = & \frac{i}{\eps} \sum_k
\left[ {\cal V}\left(\frac{t}{\eps^2},n,k\right) \rhond(t,k,n)
- {\cal V}\left(\frac{t}{\eps^2},k,n\right) \rhond(t,n,k)  \right] \\
\non 
& + & \sum_k [ W(k,n) \rhod(t,k) - W(n,k) \rhod(t,n) ],
\end{eqnarray}
and for the coherences:
\begin{eqnarray}
\label{eqnd}
\d_t \rhond(t,n,m) 
& = & \frac1{\eps^2} \Omega^\eps(n,m) \rhond(t,n,m) \\
\non
& + & \frac{i}{\eps} {\cal V}\left(\frac{t}{\eps^2},n,m\right) 
[\rhod(t,m) - \rhod(t,n)] \\
\non
& + &  \frac{i}{\eps}  \sum_k
\left[ {\cal V}\left(\frac{t}{\eps^2},n,k\right) \rhond(t,k,m)
- {\cal V}\left(\frac{t}{\eps^2},k,m\right) \rhond(t,n,k) \right].
\end{eqnarray}
\reseteqn
As a consequence of the Hermitian properties of Eq.~\eqref{bloch} recalled in
Lemma~\ref{conservation}, Eq.~\eqref{eqd} can also be cast as
\begin{eqnarray}
\label{eqdbis}
\d_t \rhod(t,n) & = & - \frac{2}{\eps} \im \sum_k
\left[ {\cal V}\left(\frac{t}{\eps^2},n,k\right) \rhond(t,k,n) \right] \\
\non
& + & \sum_k [ W(k,n) \rhod(t,k) - W(n,k) \rhod(t,n) ].
\end{eqnarray}

\subsection{An equation for populations only}

In this section, we transform the coupled system~\eqref{eqd}-\eqref{eqnd} into
one equation governing the populations $\rhod(t,n)$ only. More precisely, we
show the following proposition.

\begin{proposition}
\label{Pr_rhodun}
Define the time dependent transition rate
\[
\Psi_\eps\left(\frac{t}{\eps^2},k,n\right)
:= 2 |V(n,k)|^2 \re \int_0^{t/\eps^2} ds \exp\left(\Omega^\eps(k,n) s\right)
\phi\left(\frac{t}{\eps^2}\right) \phi\left(\frac{t}{\eps^2}-s\right).
\]
Then, for all $T>0$, the vector $\rhod$ satisfies
\begin{equation}
\label{rhod}
\d_t \rhod(t) =  
\left( \Psi_\eps\left(\frac{t}{\eps^2}\right) + W \right)\op \rhod(t) 
+ O_{L^\infty([0,T],l^1)}\left(\eps^{1-2\mu}\right).
\end{equation}
\end{proposition}

\begin{rem}
Eq.~\eqref{rhod} is a linear Boltzmann type equation with a time dependent
transition rate. This is our first description of the behavior of the
populations at leading order in $\eps$.

We could likewise obtain approximations of $\rhod$ at \textit{each order} (in
$\eps$), thus providing a hierarchy of Boltzmann type equations for the
successive orders.
\end{rem}

\begin{rem}
Using Lemma~\ref{nonpositif} and since 
$\re \Omega^\eps(n,m) \leq - \eps^{\mu} \gamma$, uniformly in $n$ and $m$, 
$n\neq m$, operator $\Psi_\eps$ is \textit{a priori} of order $\eps^{-\mu}$ 
on $l^2$ or $l^1$. More precisely we have the following estimates.
\alpheqn
\begin{eqnarray}
\label{estpsiea}
\left\| \Psi_\eps(t/\eps^2)\op \rhod \right\|_{l^1}
& \leq &
C \eps^{-\mu} \|V(n,k)\|^2_{l^\infty_n l^1_k \cap l^\infty_k l^1_n}
\|\phi\|_{L^\infty(\R)}^2 \|\rhod\|_{l^1}, \\
\label{estpsieb}
\left\| \Psi_\eps(t/\eps^2)\op \rhod \right\|_{l^2}
& \leq &
C \eps^{-\mu} \|V(n,k)\|^2_{l^\infty_n l^1_k \cap l^\infty_k l^1_n}
\|\phi\|_{L^\infty(\R)}^2 \|\rhod\|_{l^2}.
\end{eqnarray}
\reseteqn
\end{rem}

\begin{proof}

\noi
Proposition~\ref{Pr_rhodun} is proved in three steps that follow \cite{BCD}.
First, coherences $\rhond(t,n,k)$ are computed at leading order in $\eps$.
These leading order terms are expressed in terms of populations $\rhod(t,n)$
\textit{only} as stated in Lemma~\ref{Lm_coherence} below. Next, we plug this
result into Eq.~\eqref{eqdbis} governing populations. The closed equation for
populations \eqref{eqdter} which is thus obtained is a linear Boltzmann
equation with a time-delay term. To show this delay is indeed small and can be
removed to obtain the delay-free equation \eqref{rhod} we then use
calculations inspired by \cite{Ca1,Ca2,Ca3}.

\paragraph{\it First step: computation of coherences}~\\
Since the initial data for coherences is $\rhond(t=0,n,m) \equiv 0$, the
integral form for Eq.~\eqref{eqnd} reads
\begin{eqnarray}
\non
\rhond(t,n,m) & = &
i \eps \int_0^{t/\eps^2} ds \exp\left(\Omega^\eps(n,m) s\right)
{\cal V}\left(\frac{t}{\eps^2}-s,n,m\right)
[\rhod(t-\eps^2 s,m) - \rhod(t - \eps^2 s,n)]
\\
&&
\label{resorhond}
+ (A_\eps\rhond)(t,n,m),
\end{eqnarray}
where
\begin{multline*}
(A_\eps\rhond)(t,n,m) :=
i \eps \int_0^{t/\eps^2} ds \exp\left(\Omega^\eps(n,m) s\right) \times
\\
\times
\sum_k
\left[ {\cal V}\left(\frac{t}{\eps^2}-s,n,k\right) \rhond(t-\eps^2 s,k,m)
- {\cal V}\left(\frac{t}{\eps^2}-s ,k,m\right) \rhond(t-\eps^2 s,n,k)
\right].
\end{multline*}
We are only interested in the first term in the expansion of $\rhond$. Indeed,
we only use leading order terms in $\eps$. From the integral
equation~(\ref{resorhond}), we would however be able to express $\rhond$ as a
{\em complete} expansion in powers of $\eps$, in terms of $\rhod$. The
following lemma states that the remainder of the expansion is indeed small.

\begin{lemma}
\label{Lm_coherence}
Set
\[
\rhond^{(0)}(t,n,m)
:=
i \eps \int_0^{t/\eps^2} ds \exp\left(\Omega^\eps(n,m) s\right)
{\cal V}\left(\frac{t}{\eps^2}-s,n,m\right)
[\rhod(t-\eps^2 s,m) - \rhod(t - \eps^2 s,n)].
\]
Then, for any given time $T\geq 0$, there exists a constant $C$, that does not
depend on $\eps$, such that we have the estimates
\begin{eqnarray*}
\Big\| \rhond - \rhond^{(0)} \Big\|_{L^\infty([0,T],l^1)}
& \leq & C \eps^{2(1-\mu)}, \\
\Big\| \rhond \Big\|_{L^\infty([0,T],l^1)} & \leq & C \eps^{1-\mu}.
\end{eqnarray*}
\end{lemma}

\begin{proof} (Lemma~\ref{Lm_coherence}).~\\
Let $T \geq 0$ be given. Because 
$\re \Omega^\eps(n,m) \leq - \eps^\mu \gamma < 0$, uniformly in $n$ and $m$, 
$n\neq m$, we have the estimate
\begin{eqnarray*}
\|A_\eps \rhond\|_{L^\infty([0,T],l^1)}
& \leq &
2 \eps
\left\|
\int_0^{+\infty}
ds \left| \exp\left(\Omega^\eps(n,m) s\right) \right|
\right\|_{l^\infty_{n,m}}
\|{\cal V}\|_{L^\infty(\R^+,l^\infty_n l^1_m \cap l^\infty_m l^1_n)}
\|\rhond\|_{L^\infty([0,T],l^1)}
\\
& \leq & C \eps^{1-\mu} \|\rhond\|_{L^\infty([0,T],l^1)}.
\end{eqnarray*}
According to the definition of $\rhond^{(0)}(t,n,m)$, Eq.~\eqref{resorhond}
\[
\rhond(t,n,m) = \rhond^{(0)}(t,n,m) + (A_\eps\rhond)(t,n,m),
\]
and we can estimate the difference
\begin{eqnarray*}
\|\rhond- \rhond^{(0)}\|_{L^\infty([0,T],l^1)}
& \leq & C \eps^{1-\mu} \|\rhond\|_{L^\infty([0,T],l^1)}
\\
& \leq & C \eps^{1-\mu}\; \|\rhond-\rhond^{(0)}\|_{L^\infty([0,T],l^1)}
+ C \eps^{1-\mu} \|\rhond^{(0)}\|_{L^\infty([0,T],l^1)}.
\end{eqnarray*}
Hence, if $\eps$ is small enough, we have
\begin{eqnarray*}
\|\rhond - \rhond^{(0)}\|_{L^\infty([0,T],l^1)}
& \leq &
C \eps^{1-\mu} \|\rhond^{(0)}\|_{L^\infty([0,T],l^1)}.
\end{eqnarray*}
According to the definition of $\rhond^{(0)}(t,n,m)$, and using the same type
of estimates,
\begin{eqnarray*}
\|\rhond^{(0)}\|_{L^\infty([0,T],l^1)}
&\leq &
C \eps^{1-\mu}
\|{\cal V}\|_{L^\infty(\R; l^\infty_n l^1_k \cap l^\infty_k l^1_n)}
\|\rhod\|_{L^\infty([0,T],l^1)}.
\end{eqnarray*}
Now the crucial estimate stems from the trace conservation
property~\eqref{lindblad} which also reads
$\|\rhod\|_{L^\infty([0,T],l^1)} = \|\rhod(t=0)\|_{l^1}$. This $l^1$-estimate 
on the diagonal terms of the density matrix is therefore sufficient to control 
all the off-diagonal terms in turn, and Lemma~\ref{Lm_coherence} follows.
\end{proof}

\paragraph{\it Second step: the time-delayed differential equation for the 
populations}~\\
Lemma~\ref{Lm_coherence} together with Eq.~\eqref{eqdbis} governing $\rhod$
implies that
\begin{eqnarray}
\label{eqdter}
\d_t \rhod(t,n) & = & \sum_k [ W(k,n) \rhod(t,k) - W(n,k) \rhod(t,n) ] \\
\non
& + & \sum_k \int_0^{t/\eps^2} ds [ \rhod(t-\eps^2 s,k) - \rhod(t-\eps^2 s,n)] 
\times \\
\non
&&
\qquad
\times 2 \re
\left\{ \exp(\Omega^\eps(k,n) s) {\cal V}\left(\frac{t}{\eps^2},n,k\right)
{\cal V}\left(\frac{t}{\eps^2}-s,k,n\right)
\right\}
\\
\non
&& +  O_{L^\infty([0,T],l^1)}\left(\eps^{1-2\mu}\right).
\end{eqnarray}

\paragraph{\it Third step: convergence to a delay-free equation}~\\
From now on we will only deal with Boltzmann type equations and use 
extensively the shorter expressions defined in Notation~\ref{Not_op}. Hence 
Eq.~\eqref{rhod} can be cast as
\[
\d_t \rhod(t,n) = \left(\Psi_\eps(t/\eps^2)\op \rhod\right)(t,n)
+ (W\op \rhod)(t,n) + O_{L^\infty([0,T],l^1)}\left(\eps^{1-2\mu}\right).
\]
Moreover, if $T \geq 0$ is given, the delayed terms $\rhod(t-\eps^2 s)$ in
Eq.~\eqref{eqdter} read
\[
\rhod(t-\eps^2s,n) = \rhod(t,n)
+ O \left(\eps^2 s \|\d_t\rhod(\cdot,n)\|_{L^\infty([0,T])}\right).
\]
Thus, Eq.~\eqref{eqdter} yields
\begin{equation}
\label{eqdfour}
\d_t \rhod(t,n) = \left(\Psi_\eps (t/\eps^2)\op \rhod\right)(t,n)
+ (W\op \rhod)(t,n) + O_{L^\infty([0,T],l^1)}\left(\eps^{1-2\mu}\right)
+ r_\eps(t,n),
\end{equation}
where the remainder $r_\eps$ can be estimated by
\begin{align*}
\|r_\eps\|_{L^\infty([0,T],l^1)} & \leq C \|\d_t \rhod\|_{L^\infty([0,T],l^1)}
\left\|
\eps^2 \disp \int_0^{+\infty}
ds \; s \exp(\Omega^\eps(n,m) s) \right\|_{l^\infty_{n,m}} \hspace*{-5cm}\\
& \leq
C \eps^{2-2\mu} \|\d_t \rhod\|_{L^\infty([0,T],l^1)}
\\
& \leq
C \eps^{2-2\mu} \eps^{-\mu} \|\rhod\|_{L^\infty([0,T],l^1)}
&& \textrm{thanks to Eqs~\eqref{eqdfour} and \eqref{estpsiea}}
\\
& \leq
C \eps^{2-3\mu}
&& \textrm{thanks to Eq.~\eqref{lindblad}},
\end{align*}
for some constant $C$, that does not depend on $\eps$, if $\eps$ is small
enough. Including this new estimate in Eq.~\eqref{eqdfour}, we have
\[
\d_t \rhod(t,n) = \left((\Psi_\eps)\op \rhod\right)(t,n) + (W\op \rhod)(t,n) 
+ O_{L^\infty([0,T],l^1)}\left(\eps^{1-2\mu}+\eps^{2-3\mu}\right).
\]
We recall that $\mu<1/2$ therefore $1-2\mu<2-3\mu$ and 
Proposition~\ref{Pr_rhodun} follows.
\end{proof}

\section{Time-averaging of transition rates}
\label{Sec_ave}

Proposition~\ref{Pr_rhodun} reduces the problem to the asymptotic analysis of
a leading order equation, namely Eq.~\eqref{rhod}. As $\eps$ goes to zero, the 
rapid oscillations of the time-dependent coefficient $\Psi_\eps(t/\eps^2)$ are 
naturally smeared out, so that only the time average of $\Psi_\eps$ plays a 
significant r\^ole. The corresponding rigorous mathematical statement is 
proven here, as in \cite{BCD}, upon using averaging techniques (see 
\textit{e.g.} \cite{SV}). This leads to substitute Eq.~\eqref{rhod} by a new
Boltzmann equation, with time-independent rates.

The result of these averaging procedures strongly depends on the specific form 
of the wave. In the present article, we need explicit forms for the averaged 
transition rates and therefore restrict to the case of a quasiperiodic wave, 
as described in Section~\ref{Sec_fonctionnel}. The small divisor problems that 
stem from this analysis are handled assuming strong decay properties on the
Fourier coefficients, namely Hypothesis~\ref{smooth}.

The averaged transitions rates we obtain in this section still depend on
$\eps$. Besides we do not discriminate dominant and negligible (with respect
to $\eps$) contributions. We postpone this sorting out to 
Section~\ref{Sec_dom} where negligible contributions will be pointed out as 
non-resonant contributions. \\

Using the explicit value of the wave $\phi(t)$ given by Eq.~\eqref{Fourier},
we compute from Eq.~\eqref{rhod} the time-dependent transition rate
$\Psi_\eps$:
\begin{multline}
\label{psiexpli}
\Psi_\eps\left(\frac{t}{\eps^2},k,n\right)
= 2 |V(n,k)|^2 \re \sum_{\alpha,\beta \in \Z^r} \phi_\alpha \phi_\beta
\exp\left(i(\alpha+\beta)\cdot\omega \frac{t}{\eps^2}\right) \times
\\
\times \frac
{1-\exp\left([-\eps^\mu\gamma(k,n)
-i(\omega(k,n)+\beta\cdot\omega+\eps^p\delta(k,n))] t/\eps^2\right)}
{\eps^\mu\gamma(k,n)+i(\omega(k,n)+\beta\cdot\omega+\eps^p\delta(k,n))}.
\end{multline}
The goal of this section is to prove the following result.

\begin{proposition}
\label{Pr_relax}
{\bf Define the averaged transition rate
\[
\Psiave(k,n) := \lim_{s\to+\infty} \frac{1}{s} \int_0^s ds
\Psi_\eps(s).
\]
Its explicit value is } 
\begin{equation}
\label{psidomeps}
\Psiave(k,n) = 2 |V(n,k)|^2 \sum_{\beta \in \Z^r} 
\frac{\eps^\mu \gamma(k,n)}
{\eps^{2\mu}\gamma(k,n)^2 + |\omega(n,k)+\beta\cdot\omega+\eps^p\delta(k,n)|^2}
|\phi_\beta|^2.
\end{equation}
Define $\rhoun$ as the vector solution to
\begin{equation}
\label{rhod1}
\d_t \rhoun = ( \Psiave + W )\op \rhoun
\end{equation}
with initial data $\rhoun(0,n)=\rhod(0,n)$. Under Hypotheses~\ref{dioph} and
\ref{smooth}, for all $T>0$, there exists $C>0$ such that
\[
\|\rhod-\rhoun \|_{L^\infty([0,T],l^2)} \leq C \eps^{1-2\mu}.
\]
\end{proposition}

\begin{rem}
In the transition rate $\Psiave(k,n)$ there are still contributions of 
non-resonant waves, \textit{i.e.} thoses such that 
$\omega(n,k)+ \beta\cdot\omega\neq 0$. They are removed in the next section. 
However, we gained in defining a time-indepedent transition rate, in that we 
filtered oscillations.
\end{rem}

The simple but crucial remark that leads from Proposition~\ref{Pr_rhodun} to 
Proposition~\ref{Pr_relax} is the fact that
$\Psiave \in l^\infty_n l^1_k \cap l^\infty_k l^1_n$ and its entries are 
non-negative. We can therefore apply Lemma~\ref{nonpositif} directly.

\begin{lemma}
\label{negaga}
For any $\eps>0$, operator $\Psiave\op$ is a bounded \textit{non-positive} 
operator on the Hilbert space $l^2$. In particular, the exponential 
$\exp(t \Psiave\op)$ is well defined as an operator on $l^2$ for $t\geq 0$, 
and its norm is less than $1$, for all $t \geq 0$.
\end{lemma}

(Proof: see Appendix~\ref{Sec_relax}, Lemma~\ref{nonpositif})

\begin{rem}
This property proves to be crucial in the sequel, and therefore the asymptotic
result in Proposition~\ref{Pr_relax} and Theorem~\ref{Th_relax} may 
not hold when the transition rate has no sign. This is why we restrict the
analysis to the case of a quasiperiodic wave when relaxations tend to zero
with $\eps$. The asymptotic result may also be false if applied to the
time-dependent transition rate $\Psi_\eps(t/\eps^2)$, which clearly changes
signs.
\end{rem}

\begin{proof}
We follow \cite{BCD} to prove Proposition~\ref{Pr_relax}.

\paragraph{\it First step: splitting of $\Psi_\eps$}~\\
Proposition~\ref{Pr_relax} amounts to estimating the difference
\[
\Delta(t) := \rhod(t) - \rhoun(t).
\]
To this aim, we first give shorter forms for the equations governing
$\rhod$ and $\rhoun$ respectively. Namely Eq.~\eqref{rhod} can be cast as
\[
\d_t\rhod(t)
= \left( \Psi_\eps\left(\frac{t}{\eps^2}\right)\op \rhod \right)(t)
+ ( W\op \rhod )(t)
+ O\left(\eps^{1-2\mu}\right).
\]
where the transition rate $\Psi_\eps(t/\eps^2)$ is given by
Eq.~\eqref{psiexpli}, and Eq.~\eqref{rhod1} also reads
\[
\d_t\rhoun(t) = \left(\Psiave\op \rhoun\right)(t) 
+ \left(W\op \rhoun\right)(t).
\]
Hence the difference $\Delta(t)$ satisfies the equation
\begin{eqnarray}
\label{eqde}
\d_t \Delta(t) = (\Psiave\op \Delta)(t)
+ \left(\Psiosc\left(\frac{t}{\eps^2}\right)\op \rhod\right)(t)
+ (W\op \Delta)(t) + O(\eps^{1-2\mu}),
\end{eqnarray}
where
\[
\Psiosc \left(\frac{t}{\eps^2},k,n\right)
:=
\Psi_\eps\left(\frac{t}{\eps^2},k,n\right) - \Psiave(k,n)
\]
contains the oscillatory contribution to the transition rate, which we want to
prove to be negligible. Gathering the terms for which $\alpha+\beta=0$, this
contribution is equal to
\begin{eqnarray}
\non
\Psiosc \left(\frac{t}{\eps^2},k,n\right)
& = & 2 |V(n,k)|^2 \re \Bigg( 
\begin{aligned}[t]
& - \sum_{\beta \in \Z^r}
\frac{|\phi_\beta|^2}{\eps^\mu \gamma(k,n) 
+ i(\omega(k,n)+\beta\cdot\omega+\eps^p\delta(k,n))} \times \\
& \times
\exp([-\eps^\mu\gamma(k,n) -i(\omega(k,n)+\beta\cdot\omega+\eps^p\delta(k,n))]
\frac{t}{\eps^2}) 
\end{aligned}
\\
&&
\label{defpsieps}
+ \sum_{\alpha\neq-\beta \in \Z^r}
\begin{aligned}[t]
& \frac{\phi_\alpha \phi_\beta 
\exp(i(\alpha+\beta)\cdot\omega t/\eps^2)}
{[\eps^\mu\gamma(k,n)+i(\omega(k,n)+\beta\cdot\omega+\eps^p \delta(k,n))]}
\times \\
&
\times
[1 - \exp([-\eps^\mu\gamma(k,n)
-i(\omega(k,n)+\beta\cdot\omega+\eps^p\delta(k,n))]
\frac{t}{\eps^2})]
\Bigg).
\end{aligned}
\end{eqnarray}
This expression carries ``time-oscillations'' (at frequency $\eps^{-2+\mu}$
at least), which kill the possibly diverging factors $\eps^{-\mu}$ (due to 
the denominators), and make them of size $\eps^{2-2\mu}$.

\paragraph{\it Second step: preliminary bounds}~\\
Since $V\in l^\infty_n l^1_k \cap l^\infty_k l^1_n$ and
$\sum_\beta |\phi_\beta|^2 < \infty$, we first find that
\[
\begin{array}{ll}
\disp \sup_n \sum_{k} |\Psiave(k,n)| \leq C \eps^{-\mu}, \quad &
\disp \sup_n \sum_{k} |\Psiosc(t/\eps^2,k,n)| \leq C \eps^{-\mu},
\\
\disp \sup_k \sum_{n} |\Psiave(k,n)| \leq C \eps^{-\mu}, \quad &
\disp \sup_k \sum_{n} |\Psiosc(t/\eps^2,k,n)| \leq C \eps^{-\mu},
\end{array}
\]
for some $C>0$ that does not depend on $t$ and $\eps$. Lemma~\ref{nonpositif}
yields the operator estimates
\begin{equation}
\label{bdpsidomop}
\left\|\Psiave\op u \right\|_{l^2} \leq C \eps^{-\mu} \|u\|_{l^2},
\quad
\left\|{\Psiosc}\op u \right\|_{l^2} \leq C \eps^{-\mu} \|u\|_{l^2}.
\end{equation}
Besides, we have the upper bound 
\begin{eqnarray}
\label{ass}
\sup_{0 \leq t \leq T}
\left\| \int_0^{t/\eps^2} ds \Psiosc(s,k,n)
\right\|_{l^\infty_n l^1_k \cap l^\infty_k l^1_n}
\leq C \eps^{-2\mu}.
\end{eqnarray}
Now, according to Eq.~\eqref{defpsieps}, $\Psiosc(t,k,n)$ is a sum of two
different terms. We use the decay assumptions
$V \in l^\infty_n l^1_k \cap l^\infty_k l^1_n$ and $\phi_\alpha\in l^2$ to
estimate the contribution of the first term to the integral by
\[
C \left\|
\sum_\beta |V(n,k)|^2 \frac{|\phi_\beta|^2}{|\eps^\mu\gamma(k,n)
+i(\omega(k,n)+\beta\cdot\omega+\eps^p\delta(k,n)|^2}
\right\|_{l^\infty_n l^1_k \cap l^\infty_k l^1_n}
\leq C \eps^{-2 \mu}.
\]
The second contribution is estimated by
\begin{eqnarray*}
&&
C \left\|
\sum_{\alpha+\beta\neq 0}
\frac{|V(n,k)|^2 |\phi_\alpha| |\phi_\beta|}
{|\eps^\mu \gamma(k,n)
+ i (\omega(k,n)+\beta\cdot\omega+\eps^p \delta(k,n))|}
\cdot \frac{1}{|(\alpha+\beta)\cdot\omega|}
\right\|_{l^\infty_n l^1_k \cap l^\infty_k l^1_n}
\\
&&
\qquad
+ C \left\|
\sum_{\alpha+\beta\neq 0}
\frac{|V(n,k)|^2 |\phi_\alpha| |\phi_\beta|}
{|\eps^\mu \gamma(k,n)+i ((\alpha+\beta)\cdot\omega+\eps^p \delta(k,n))|^2}
\right\|_{l^\infty_n l^1_k \cap l^\infty_k l^1_n}
\\
&&
\leq C \eps^{-\mu} \sum_{\alpha,\beta} |\phi_\alpha| |\phi_\beta| 
|\alpha+\beta|^{r-1+\eta}
+ C \eps^{-2\mu} \sum_{\alpha, \beta} |\phi_\alpha| |\phi_\beta|
\\
&& \leq C \eps^{-2 \mu},
\end{eqnarray*}
thanks to the Diophantine estimate (Hypothesis~\ref{dioph}), together with
Hypothesis~\ref{smooth}. This yields inequality \eqref{ass}.

\paragraph{\it Third step: integral form of the equations}~\\
Since $\Delta(0)=0$, the integral form for Eq.~\eqref{eqde} governing 
$\Delta(t)$ reads
\[
\Delta(t) := \Delta^{(1)}(t)+\Delta^{(2)}(t),
\]
where
\begin{eqnarray*}
\Delta^{(1)}(t) & = &
\int_0^t ds \exp([t-s] \Psiave)
\Psiosc\left(\frac{s}{\eps^2}\right)\op \rhod(s), \\
\Delta^{(2)}(t) & = &
\int_0^t ds \exp([t-s] \Psiave) 
\left( (W\op \Delta)(s) + O(\eps^{1-2\mu}) \right).
\end{eqnarray*}

\paragraph{\it Fourth step: estimating $\Delta^{(1)}(t)$ and 
$\Delta^{(2)}(t)$}~\\
Here Lemma~\ref{negaga} proves crucial, in that we use the exponential of the
bounded operator $\Psiave$. This together with Lemma~\ref{nonpositif} applied
to operator $W$ yields an estimate for $\Delta^{(2)}$:
\begin{eqnarray}
\label{gron}
\|\Delta^{(2)}(t)\|_{l^2}
\leq
C \left( \eps^{1-2\mu} + \int_0^t ds \|\Delta(s)\|_{l^2} \right).
\end{eqnarray}
On the other hand, to take advantage of the time oscillations of the operator
$\Psiosc(t/\eps^2)$, we carry out a natural integration by parts in the
expression for $\Delta^{(1)}$:
\begin{eqnarray*}
\Delta^{(1)}(t) & = &
\eps^2 \left( \int_0^{t/\eps^2} du \Psiosc(u) \right)\op \rhod(t)
\\
& + &
\eps^2 \int_0^t ds \exp([t-s] \Psiave) \Psiave
\left( \int_0^{s/\eps^2} du \Psiosc(u) \right)\op \rhod(s)
\\
& - &
\eps^{2} \int_0^t ds \exp([t-s] \Psiave)
\left( \int_0^{s/\eps^2} du \Psiosc(u) \right) \times
\\
&&
\qquad\qquad
\times \left( \Psiave + \Psiosc\left(\frac{s}{\eps^2}\right)
+ W + O(\eps^{1-2\mu}) \right)\op \rhod(s),
\end{eqnarray*}
where we have used Eq.~\eqref{rhod} to express $\d_t \rhod(s)$.
Estimates~\eqref{bdpsidomop} on the operators $\Psiave$ and $\Psiosc(t)$
together with Lemma~\ref{nonpositif} and Lemma~\ref{negaga} (non-positiveness
of $\Psiave$) lead to
\[
\|\Delta^{(1)}\|_{L^\infty([0,T],l^2)}
\leq C  \eps^{2-\mu} \sup_{0\leq t \leq T}
\left\| \int_0^{t/\eps^2} ds \Psiosc(s) \right\|_{{\cal L}(l^2)}
\|\rhod\|_{L^\infty([0,T],l^2)}.
\]
Besides we have
$\|\rhod\|_{L^\infty([0,T],l^2)} \leq \|\rhod\|_{L^\infty([0,T],l^1)} \leq C$,
and it follows that, for $\eps$ small enough, we have
\begin{eqnarray*}
\|\Delta^{(1)}\|_{L^\infty([0,T],l^2)}
& \leq &
C  \eps^{2-\mu} \sup_{0\leq t \leq T}
\left\| \int_0^{t/\eps^2} ds \Psiosc(s) \right\|_{{\cal L}(l^2)}
\\
& \leq &
C  \eps^{2-\mu} \sup_{0\leq t \leq T}
\left\| \int_0^{t/\eps^2} ds \Psiosc(s,k,n)
\right\|_{l^\infty_n l^1_k \cap l^\infty_k l^1_n}
\\
& \leq & C \eps^{2-3\mu}.
\end{eqnarray*}
This, together with estimate~\eqref{gron}, and Gronwall lemma, yields
\[
\|\Delta(t)\|_{L^\infty([0,T],l^2)} \leq C \eps^{1-2\mu}.
\]
and Proposition~\ref{Pr_relax} is proved.
\end{proof}

\section{Keeping only resonant contributions in transition\break rates:
the main theorem}
\label{Sec_dom}

In this section we prove that the non-resonant contributions, which correspond
to the triples $(n,k,\beta)$ such that $\omega(n,k)+\beta\cdot\omega\neq 0$ in 
the transition rate~\eqref{psidomeps}, are negligible in the limit $\eps\to0$. 
We therefore replace the transition rate $\Psiave$ by a purely resonant
transition rate $\Psidom$, which however still depends on $\eps$. To get rid
of this last dependence, we will have to specify the value of $p/\mu$, which
is the goal of the next section.

Still due to small denominator problems, we need to reinforce the decay
assumptions on the coefficients and assume Hypothesis~\ref{extrasmooth} and
\ref{far} hold.

\begin{theorem}
\label{Th_relax}
Define the transition rate
\begin{equation}
\label{psidom}
\Psidom(k,n) 
:= 2 |V(n,k)|^2 \frac{\eps^\mu \gamma(k,n)}{\eps^{2\mu}\gamma(k,n)^2 
+ \eps^{2 p} \delta(k,n)^2}
\sum_{ \beta \in \Z^r \suchthat \omega(k,n)+\beta\cdot\omega=0} |\phi_\beta|^2.
\end{equation}
Let also $\rhodeux$ be solution to
\begin{equation}
\label{rhod2}
\d_t \rhodeux = \left( \Psidom + W \right)\op \rhodeux
\end{equation}
with initial data $\rhodeux(0,n)=\rhod(0,n)$. We assume that $\mu < 1/2$.
Then, under Hypotheses~\ref{dioph}, \ref{speed}, \ref{extrasmooth} and 
\ref{far}, for all $T>0$, there exists $C>0$ such that
\[
\|\rhod-\rhodeux\|_{L^\infty([0,T],l^2)} \leq C (\eps^\mu + \eps^{1-2\mu}).
\]
\end{theorem}

\begin{rem}
When $\mu=0$, Theorem~\ref{Th_relax} does not give a good approximation of 
$\rhod$, and we do not have a better description than the one from 
Proposition~\ref{Pr_relax}, where all (resonant and non-resonant) frequencies 
have a contribution to the transition rates.
\end{rem}

The proof of Theorem~\ref{Th_relax} follows closely that of
Proposition~\ref{Pr_relax}. However, we notice that the resonant values that
play a r\^ole here satisfy $\omega(n,k)+\beta\cdot\omega=0$ and not
$\omega(n,k)+\beta\cdot\omega+\eps^p \delta(n,k)=0$. Therefore we need to 
understand the effect of a perturbation on Diophantine estimates. 
Lemma~\ref{diophstable} below answers this problem.

\subsection{Perturbed Diophantine estimates}

\begin{lemma}
\label{diophstable}
If $\omega$ and $\omega(n,k)$ satisfy the Diophantine condition~\eqref{diopha} 
and Hypothesis~\ref{speed} with constants $\eta$ and $C_\eta$, then the 
following assertion holds.

If $(n,k,\beta) \in \N\times\N\times\Z^d$ satisfies
\[
|\beta \cdot \omega + \omega(n,k)+\eps^p \delta(n,k)|
\leq
\frac12 \frac{C_\eta}{(1+|\beta|)^{r-1+\eta} (1+n)^{1+\eta} (1+k)^{1+\eta}},
\]
then
\[
(1+|\beta|)^{r-1+\eta} (1+n)^{1+\eta} (1+k)^{1+\eta} \geq
\frac{C_\eta \eps^{-p}}{2|\delta|_{l^\infty}}.
\]
\end{lemma}

\begin{rem}
The Diophantine condition~\eqref{diopha} is not stable with respect to small 
perturbations: coefficients $\omega(n,k)+\eps^p\delta(n,k)$, that can be
arbitrarily close to $\omega(n,k)$, are capable of violating the Diophantine
condition~\eqref{diopha}. Indeed, for $\eps$ small, the condition
\[
|\beta \cdot \omega + \omega(n,k)+\eps^p \delta(n,k)|
\geq
\frac12 \frac{C_\eta}{(1+|\beta|)^{r-1+\eta} (1+n)^{1+\eta} (1+k)^{1+\eta}}
\]
does therefore not necessarily hold, only assuming that the left hand-side is
non-zero. Nevertheless, Lemma~\ref{diophstable} claims that this condition may
only be violated for values of the triple $(n,k,\beta)$ which are very large 
when $\eps\to0$. See \cite{CCCFLLT} for a similar argument.
\end{rem}

\begin{proof}

\noi
Set $K=C_\eta/2$ and take $(n,k,\beta)$ such that
\[
|\beta \cdot \omega + \omega(n,k)+\eps^p \delta(n,k)|
\leq
\frac{K}{(1+|\beta|)^{r-1+\eta} (1+n)^{1+\eta} (1+k)^{1+\eta}}.
\]
Then
\[
|\beta \cdot \omega + \omega(n,k)| - \eps^p |\delta(n,k)|
\leq
\frac{K}{(1+|\beta|)^{r-1+\eta} (1+n)^{1+\eta} (1+k)^{1+\eta}},
\]
and according to condition~\eqref{diopha} (or Hypothesis~\ref{speed},
when $\beta=0$)
\[
\frac{2 K}{(1+|\beta|)^{r-1+\eta} (1+n)^{1+\eta} (1+k)^{1+\eta}}
- \eps^p |\delta(n,k)|
\leq
\frac{K}{(1+|\beta|)^{r-1+\eta} (1+n)^{1+\eta} (1+k)^{1+\eta}}.
\]
Hence
\[
\frac{K}{(1+|\beta|)^{r-1+\eta} (1+n)^{1+\eta} (1+k)^{1+\eta}}
\leq
\eps^p |\delta(n,k)|,
\]
which ends the proof of Lemma~\ref{diophstable}.
\end{proof}

\subsection{Proof of the main theorem}

\paragraph{\it First step: Integral formulation}~\\
To deduce Theorem~\ref{Th_relax} from Proposition~\ref{Pr_relax} only amounts
to estimate the difference
\[
\Delta(t) = \rhodeux(t)-\rhoun(t),
\]
where we use once again the notation $\Delta$. Now $\rhodeux(t)$ and 
$\rhoun(t)$ are respectively solution to
\[
\d_t\rhodeux(t) 
= \left(\Psidom\op \rhodeux\right)(t) + \left(W\op \rhodeux\right)(t),
\]
and
\[
\d_t\rhoun(t) 
= \left(\Psiave\op \rhoun\right)(t) + \left(W\op \rhoun\right)(t).
\]
Hence
\begin{equation}
\label{deencore}
\d_t\Delta(t) = (\Psiave\op \Delta)(t) + \left(\Psineg\op \rhodeux\right)(t) 
+ (W\op \Delta)(t).
\end{equation}
where
\[
\Psineg(k,n) = \Psidom(k,n) - \Psiave(k,n)
\]
contains the contributions to the transition rate, that we want to prove to be
negligible. Since $\Delta(0)=0$ the integral form for \eqref{deencore} reads
\[
\Delta(t)
=  \int_0^t ds \exp([t-s] \Psiave) \left(\Psineg\op \rhodeux\right)(s)
+  \int_0^t ds \exp([t-s] \Psiave) (W\op \Delta)(s).
\]

\paragraph{\it Second step: Estimating $\Psineg$}~\\
In view of Eqs~\eqref{psidomeps} and \eqref{psidom}, we have
\[
\|\Psineg\|_{l^\infty_n l^1_k\cap l^\infty_k l^1_n}
=
\sup_n \sum_{k,\beta \suchthat \omega(n,k)+\beta\cdot\omega\neq 0}
\frac{2|V(n,k)|^2 \eps^\mu \gamma(k,n) |\phi_\beta|^2}
{\eps^{2\mu}\gamma(k,n)^2
+ |\omega(k,n)+\beta\cdot\omega+\eps^p\delta(k,n)|^2}.
\]
We split this expression into two contributions according to the fact that
\[
|\beta \cdot \omega + \omega(n,k) + \eps^p\delta(n,k)|
\geq
\frac12 \frac{C_\eta}{(1+|\beta|)^{r-1+\eta} (1+n)^{1+\eta} (1+k)^{1+\eta}},
\]
or not. Using Lemma~\ref{diophstable} for the second contribution, we obtain
\begin{multline*}
\|\Psineg\|_{l^\infty_n l^1_k\cap l^\infty_k l^1_n} \\
\begin{aligned}
& \leq
\sup_n \bigg\{ \sum_{k,\beta \suchthat \omega(n,k)+\beta\cdot\omega\neq 0}
\frac{4 |V(n,k)|^2 \eps^\mu \gamma(k,n) |\phi_\beta|^2} {C_\eta}
(1+|\beta|)^{r-1+\eta} (1+n)^{1+\eta} (1+k)^{1+\eta}
\\
& \phantom{\leq \sup_n ~}
+ \sum_{k,\beta \suchthat \omega(n,k)+\beta\cdot\omega\neq 0}
\un \Big[ (1+|\beta|)^{r-1+\eta} (1+n)^{1+\eta} (1+k)^{1+\eta} \geq C
\eps^{-p} \Big] \frac{2|V(n,k)|^2 \eps^{-\mu} |\phi_\beta|^2}{\gamma}
\bigg\}.
\end{aligned}
\end{multline*}
The first sum is estimated using Hypotheses~\ref{extrasmooth}, \ref{speed} and 
\ref{far}. The second term is first multiplied and divided by the quantity
$[(1+|\beta|)^{r-1+\eta} (1+n)^{1+\eta} (1+k)^{1+\eta}]^{N_\eta}$. 
Therefore we get
\begin{eqnarray*}
\|\Psineg\|_{l^\infty_n l^1_k\cap l^\infty_k l^1_n}
& \leq &
C \eps^\mu \\
& + &
C \eps^{N_\eta p -\mu} \sum_{n,k,\beta}
\left((1+|\beta|)^{r-1+\eta} (1+n)^{1+\eta} (1+k)^{1+\eta} \right)^{N_\eta}
|V(n,k)|^2 |\phi_\beta|^2.
\end{eqnarray*}
Since we assumed that $N_\eta>2\mu/p$, we finally have
\[
\|\Psineg\|_{l^\infty_n l^1_k\cap l^\infty_k l^1_n} \leq C \eps^\mu.
\]

\paragraph{\it Third step: Conclusion}~\\
Using Lemma~\ref{negaga} (non-positiveness of $\Psiave$) leads to
\begin{eqnarray*}
\|\Delta(t)\|_{l^2}
&\leq & C \|\Psineg\|_{{\cal L}(l^2)} \| \rhodeux \|_{L^\infty([0,T],l^2)}
+ C \int_0^t ds \|\Delta(s)\|_{l^2}
\\
&\leq &
C \|\Psineg(n,k)\|_{l^\infty_n l^1_k\cap l^\infty_k l^1_n}
\left[ \|\Delta\|_{L^\infty([0,T],l^2)} + \|\rhoun\|_{L^\infty([0,T],l^2)}
\right]
+ C \int_0^t ds \|\Delta(s)\|_{l^2}.
\end{eqnarray*}
A by-product of Proposition~\ref{Pr_relax} is that the quantity
$\|\rhoun\|_{L^\infty([0,T],l^2)}$ can be estimated by a constant. Therefore
\begin{eqnarray*}
\|\Delta(t)\|_{l^2}
& \leq & C \|\Psineg(n,k)\|_{l^\infty_n l^1_k\cap l^\infty_k l^1_n}
(1 + \| \Delta\|_{L^\infty([0,T],l^2)})
+ C \int_0^t ds \|\Delta(s)\|_{l^2}
\\
& \leq & C \eps^\mu (1 + \| \Delta\|_{L^\infty([0,T],l^2)})
+ C \int_0^t ds \|\Delta(s)\|_{l^2}.
\end{eqnarray*}
Thanks to the Gronwall lemma
\[
\|\Delta\|_{L^\infty([0,T],l^2)} \leq C \eps^\mu,
\]
which proves Theorem~\ref{Th_relax}.
\qed

\section{Time-layers and equilibrium states in the $\eps \to 0$ limit}
\label{Sec_lim}

In the previous sections we have derived the rate equation \eqref{rhod2}, with
the transition rates 
\[
\Wmod(n,m) = \Psidom(n,m) + W(n,m),
\]
which may be considered as modified rates, \textit{via} the interaction with 
the wave. We now turn to the study of the dynamics of the solution to Eq. 
\eqref{rhod2}.

For a fixed $\eps$, this dynamics is described in Appendix~\ref{Sec_relax}. 
However, in most cases, the coefficients prove to be singular in $\eps$. 
Therefore the time evolution of the solution obeys two different regimes as 
$\eps$ goes to zero: first a time-layer, and then relaxation to an 
"equilibrium" (if it exists, see Appendix \ref{Sec_relax}). The duration of 
the time-layer is always less than $O(1)$, thus the Bloch equation 
\eqref{bloch} and the rate equation \eqref{rhod2} behave in the same way on 
this time-range (as Theorem~\ref{Th_relax} asserts). On the contrary, the
relaxation towards an equilibrium state, as time goes to infinity, is
\textit{a priori} specific for the rate equation, since
Theorem~\ref{Th_relax} only applies on the fixed time-interval $[0,T]$. 

In this section, we split the rate operator $(\Wmod)\op$ into three 
contributions.
\begin{enumerate}
\item Some levels are decoupled from all other levels. The corresponding 
columns (and subsequently lines) in the matrix $(\Wmod)\op$ are identically 
zero. These levels have a constant population for all time (exactly for the 
rate equation, and at leading order for the Bloch equation). They are excluded 
from the sequel of the argument. We call $\Sigma_0$ the subspace spanned by 
these levels in $l^2$ (this subspace does not depend on $\eps$). The 
projection onto the levels which have a non trivial time evolution is denoted 
by $\Pi_0$. Therefore $\Sigma_0=\Ker\Pi_0$.
\item We split $\Psidom$ into two contributions: $\Psising$ collects the 
\emph{singular} coefficients, which go to infinity as
$\eps\to0$. Then, $\Psinonsing$ gathers the non-singular coefficients, 
which are $O(1)$. The corresponding case study according to 
the value of $\mu/p$ is performed in Section~\ref{Sec_5setting}. 
\item Finally we define the $l^2$-orthogonal projection $\Pi$ onto the
  space $l^2$-orthogonal to $\Sigma_0$ in $\Ker\Psising$. It does not depend on
$\eps$ either. The quantity $(1-\Pi)\rhod$ will be proved to vanish in time 
$O(\eps^\sigma)$, where $\sigma$ depends on the ratio $\mu/p$ (see 
Table~\ref{table1}). After this time, \textit{i.e.} past the initial
time-layer, a polarized solution $\rhoout$ persists. It is given by
the $\eps$-independent system
\[
\Pi \rhoout = \rhoout, \hspace{1cm} 
\partial_t \rhoout = \Pi \Wout\op \Pi \rhoout,
\]
where $\Wout = W + \Psi_0^{\rm nonsing}$.
\end{enumerate}

This section is organized as follows: in Section~\ref{Sec_5setting} we set some
notations and explain why the transition rates are in general splitted into
three types of terms. Their respective sizes are given by powers of
$\eps$. Under the 
assumption that the number $N$ of quantum levels is finite, a precise 
description of the time layer is given according to the value of $\mu/p$ in 
Section~\ref{Sec_5cas}. The key lemma for this analysis is given in 
Section~\ref{Sec_5lemme}. Section~\ref{Sec_5mu0} is devoted to the
particular case $\mu=0$. In Section~\ref{Sec_5infini}, we discuss the few 
cases (in terms of values of $\mu/p$) when we may conclude for an infinite 
number of levels. A rigorous framework for the justification of the 
restriction to a finite number of levels is given in Section~\ref{Sec_fini}.

\subsection{Setting for the time-layer}
\label{Sec_5setting}

We first split $\Psidom(n,m)$ into two contributions according to whether
$\delta(n,m)=0$ or not. We restrict the discussion to the case when $\mu>0$ 
since the form of the transition rates in Section~\ref{Sec_dom} does not apply 
when $\mu=0$. This case is treated separately in Section~\ref{Sec_5mu0}. To 
simplify notations, we set
\[
C(n,m) = 2 |V(n,m)|^2 
\sum_{\beta \in \Z^r \suchthat \omega(n,m)+\beta\cdot\omega=0}|\phi_\beta|^2,
\]
and $\Psidom(n,m) = A_\eps(n,m) + B_\eps(n,m)$, where
\alpheqn
\begin{eqnarray}
\label{Adom}
A_\eps(n,m) & := & C(n,m) \frac{\eps^\mu \gamma(n,m)}
{\eps^{2\mu}\gamma(n,m)^2 + \eps^{2 p} \delta(n,m)^2} \un(\delta(n,m)=0), \\
\label{Bdom}
B_\eps(n,m) & := & C(n,m) \frac{\eps^\mu \gamma(n,m)}
{\eps^{2\mu}\gamma(n,m)^2 + \eps^{2 p} \delta(n,m)^2} \un(\delta(n,m)\neq0).
\end{eqnarray}
\reseteqn
Of course Eq.~(\ref{Adom}) also reads
\[
A_\eps(n,m) = \eps^{-\mu} C(n,m) \frac1{\gamma(n,m)} \un(\delta(n,m)=0)
=: \eps^{-\mu} A(n,m).
\]
In a similar way we rewrite Eq.~(\ref{Bdom}) as
\[
B_\eps(n,m) =: \eps^{-\nu} B^\eps(n,m),
\]
where $\nu$ is chosen such that $B^\eps \to B^0$ in $l^\infty l^1$, as
$\eps\to0$, and therefore also as an operator on $l^2$. To this aim we 
consider two cases, namely $\mu\leq p$ and
$\mu\geq p$. We set
\begin{equation*}
\begin{aligned}
& \nu = \mu \textrm{ and } B^\eps(n,m) := C(n,m)
\frac{\gamma(n,m)}{\gamma(n,m)^2 + \eps^{2(p-\mu)} 
\delta(n,m)^2} \un(\delta(n,m)\neq0)
& \textrm{ if } \mu\leq p, \\
& \nu = 2p-\mu \textrm{ and } B^\eps(n,m) :=  C(n,m) 
\frac{\gamma(n,m)}{\eps^{2(\mu-p)} \gamma(n,m)^2 + \delta(n,m)^2}
\un(\delta(n,m)\neq0)
& \textrm{ if } \mu\geq p.
\end{aligned}
\end{equation*}
With the above notations, we have cast Eq.~(\ref{rhod2}) governing $y=\rhodeux$
as
\begin{eqnarray}
\label{ed}
\d_t y =( \eps^{-\mu} A + \eps^{-\nu} B^\eps + W )\op y.
\end{eqnarray}
Using the block decomposition described in Appendix~\ref{Sec_relax}, we notice 
that $\Ker(B^\eps)=\Ker(B^0)$ is constant with respect to $\eps\geq0$. This is
\textit{a priori} only true if all the minimal stable eigenspaces are 
finite-dimensional.

Now, sorting between singular and non-singular contribution yields the 
different cases in Table~\ref{table1}. In the framework of \cite{BCD} $B_\eps$ 
is identically zero. If $\mu\leq p$ then $B_\eps$ contributes to singular 
terms at the same order $\eps^{-\mu}$ as $A$. If $p<\mu<2p$, there are two 
orders of magnitude $\eps^{-\mu}$ and $\eps^{-\nu}$ in the singular term. 
Finally if $2p\leq\mu$, the effect of $B_\eps$ is of order $O(1)$ and must be 
associated to that of $W$. In other words, if $\mu<2p$, the non-zero entries 
in $A$ and $B_\eps$ contribute to $\Psising$, and $\Psinonsing=0$. If 
$\mu \geq 2p$, only the non-zero entries in $A$ contribute to $\Psising$: 
$\Psinonsing=B_\eps$, which has a non-vanishing value $\Psi_0^{\rm nonsing}$ 
only if $\mu=2p$.

\subsection{A finite dimensional lemma}
\label{Sec_5lemme}

\begin{lemma}
\label{minispectre}
Let $0\leq\nu\leq\mu$. Let $A$ and $B^\eps\in {\cal M}_N(\R)$ be symmetric 
non-positive matrices such that $B^\eps\to B^0$, and assume that 
$\Ker(B^\eps)$ is constant for $\eps\geq0$. Let $\Pi$ be the orthogonal 
projection onto $\Ker A \cap\Ker B^0$. Then there exists a constant $c>0$, 
namely
\[
c=\min_{\|x\|\leq 1, \Pi x=0} - ((A+B^0)x,x),
\]
such that any non-zero eigenvalue $\lambda^\eps$ of
$(1-\Pi) (\eps^{-\mu} A + \eps^{-\nu}B^\eps) (1-\Pi)$ satisfies
\[
\lambda^\eps \leq - c \eps^{-\nu}.
\]
\end{lemma}

\begin{proof}
Since $\mu \geq \nu$, we write
\[
(1-\Pi) (\eps^{-\mu} A + \eps^{-\nu}B) (1-\Pi)
= \eps^{-\nu} (1-\Pi) (B + \eps^{\nu -\mu}A) (1-\Pi).
\]
Let $\kappa^\eps$ be a non-zero eigenvalue of 
$(1-\Pi) (B^\eps + \eps^{\nu -\mu}A) (1-\Pi)$. Since $A$ and $B^\eps$
are non-positive matrices, we have 
$\kappa^\eps < 0$. We therefore want to prove that for the constant 
$c$ defined in the lemma, $\kappa^\eps \leq -c$. Let 
$x^\eps\in\range(1-\Pi)$ be an eigenvector associated with the eigenvalue
$\kappa^\eps$ and such that $\|x^\eps\| = 1$. We have
$(1-\Pi) (B^\eps + \eps^{\nu -\mu}A) (1-\Pi) x^\eps=\kappa^\eps x^\eps$. 
Taking the scalar product with $x^\eps$, we obtain
\[
\kappa^\eps \|x^\eps\|^2 = ((B^\eps + \eps^{\nu-\mu}A) x^\eps,x^\eps)
= (B^\eps x^\eps,x^\eps)+\eps^{\nu -\mu}(Ax^\eps,x^\eps).
\]
Since $\eps \in [0,1]$, $\nu-\mu\leq0$ and $(Ax^\eps,x^\eps) \leq 0$, we have
$\eps^{\nu -\mu}(Ax^\eps,x^\eps) \leq (Ax^\eps,x^\eps)$. A subsequence of
$x^\eps$ converges to $x^0$ and therefore
\[
\sup_{\eps\geq0} \kappa^\eps \|x^\eps\|^2 \leq ((A+B^0)x^0,x^0) 
\leq -c \|x^0\|^2.
\]
Finally, if $c$ were zero, then $((A+B^0)x,x)=0$ that is $(Ax,x)=0$ and
$(B^0x,x)=0$ and therefore $x \in \Ker A  \cap \Ker B^0$. Thus $(1-\Pi) x=0$ 
and $x=0$ which is impossible.
\end{proof}


\begin{rem}
~\\[-8mm]
\begin{itemize}
\item[(i)] The above proof is restricted to the finite dimensional case for two
reasons. First, in the infinite dimensional case the convergence of
$x^\eps$ to $x^0$ would only be weak. Second, the maximum of
$((A+B^0)x,x)$ on vectors $x$ such that $\|x\|=1$ could be zero.
\item[(ii)] To avoid such limitations we could think about replacing $B^\eps$ 
by the leading order terms (those which are not vanishing as $\eps\to0$) in 
the series expansion of $B_\eps$ (Lemma~\ref{minispectre} can clearly be 
extended to a finite number of matrices) extending the computations of
Section~\ref{Sec_5infini} to sums of several powers of $\eps$. This is not
possible since such a procedure would not ensure the non-positiveness of the
resulting operators, which is crucial in the proof.
\end{itemize}
\end{rem}

Lemma~\ref{minispectre} allows to split the solution $y$ to Eq.~(\ref{ed})
into two parts: $(1-\Pi) y$ vanishes exponentially and a solution $z$
associated with the initial data $\Pi y(0)$ survives.

\begin{Corollary}
\label{coro}
Let $0 < \nu \leq \mu$. Let $A, B^\eps \in {\cal M}_N(\R)$ be two symmetric 
non-positive matrices such that $B^\eps\to B^0$, and assume that 
$\Ker(B^\eps)$ is constant for $\eps\geq0$. Let $W \in {\cal M}_N(\R)$, and 
let $\Pi$ be the orthogonal projection onto $\Ker A\op \cap \Ker B^\eps\op$. 
If $y$ is solution to
\[
\partial_t y =(\eps^{-\mu}A+\eps^{-\nu}B^\eps + W)\op y
\]
and if $z$ is solution to 
\begin{equation}
\label{ed_z}
\partial_t z= \Pi W \op \Pi z, \qquad z(0) = \Pi y(0) ,
\end{equation}
then 
\alpheqn
\begin{eqnarray}
\label{vanish}
\|(1-\Pi)y\| & \leq & C(\eps^{\nu} + \exp(-ct\eps^{-\nu})), \\
\label{projecteur}
\|\Pi(y-z)\| & \leq & C(\eps^{\nu} + \exp(-ct\eps^{-\nu})),
\end{eqnarray}
\reseteqn
where $C$ depends on $W$, $y(0)$, $T$ and $c$ (the constant of the previous 
lemma). In particular, $\lim_{\eps\to0} (1-\Pi)y = 0$.
\end{Corollary}

\begin{proof}
\noindent
The proof is standard. We reproduce it for the sake of completeness. We first 
prove estimate \eqref{vanish}. Since $\Pi$ is the orthogonal projection onto 
$\Ker A \cap \Ker B^\eps$, we have 
$\Pi A\op  = \Pi B^\eps\op  = A\op  \Pi = B^\eps\op  \Pi = 0$, and therefore  
\begin{eqnarray*}
\partial _t (1-\Pi) y = (1-\Pi) \partial_t y 
& = & (1-\Pi) (\eps^{-\mu}A + \eps^{-\nu}B^\eps)\op y  + (1-\Pi)W\op y \\
& =  & (\eps^{-\mu}A + \eps^{-\nu}B^\eps)\op y + (1-\Pi)W\op y, \\
& = & (\eps^{-\mu}A + \eps^{-\nu}B^\eps)\op (1-\Pi)y + (1-\Pi)W\op y.
\end{eqnarray*}
The solution to this equation reads 
\begin{eqnarray*}
(1-\Pi) y(t) & = & \exp(t(\eps^{-\mu}A + \eps^{-\nu}B^\eps)\op ) (1-\Pi)y(0) \\
&& + \int_0^t ds \exp((t-s) (\eps^{-\mu}A + \eps^{-\nu}B^\eps)\op ) 
(1-\Pi)W\op y(s).
\end{eqnarray*} 
Lemma \ref{minispectre} yields a bound for the spectrum of 
$(\eps^{-\mu}A + \eps^{-\nu}B^\eps)\op$, namely 
${\rm Sp}(\eps^{-\mu}A + \eps^{-\nu}B^\eps)\op \leq -c \eps^{-\nu}$ and we 
also have $\|1-\Pi\| \leq 1$, therefore
\[
\|(1-\Pi)y \| 
\leq \exp(-ct\eps^{-\nu}) \|y(0)\| + C \int_0^t ds \exp(-(t-s)c \eps^{-\nu}).
\] 
Finally $\|(1-\Pi)y\| \leq C (\eps^{\nu} + \exp(-ct\eps^{-\nu}))$ which is 
estimate \eqref{vanish}.

We now prove estimate \eqref{projecteur}. We write 
$y-z=\Pi(y-z)+(1-\Pi)y-(1-\Pi)z$. By definition $\Pi z=z$ and we already have
estimated $(1-\Pi)y$. Thus there remains to estimate the quantity $\Pi(y-z)$, 
which is solution to 
\begin{eqnarray*}
\partial_t \Pi(y-z) & = & \partial_t \Pi y - \partial_t \Pi z
= \Pi \partial_t y - \partial_t z \\ 
& = & \Pi (\eps^{-\mu}A + \eps^{-\nu}B)\op y + \Pi W\op  y - \Pi W\op \Pi z \\ 
& = & \Pi W\op y - \Pi W\op \Pi z \\ 
& = & \Pi W\op \Pi y + \Pi W\op (1-\Pi)y - \Pi W\op \Pi z \\
& = & \Pi W\op \Pi\ \Pi(y-z)+\Pi W\op (1-\Pi)y.
\end{eqnarray*}
The solution to this equation may be estimated by
\[
\|\Pi(y-z)\| \leq \int_0^t ds 
\|\exp((t-s)\Pi W\op \Pi)\| \times \|\Pi W\op (1-\Pi)y(s)\|
\] 
because $\Pi y(0)=z(0)=\Pi z(0)$. Since $\Pi W\op \Pi \leq 0$ we finally have
\[
\|\Pi(y-z)\| \leq C (\eps^{\nu} + \exp(-ct\eps^{-\nu}))
\] 
and estimate \eqref{projecteur} follows.
\end{proof}

\begin{rem}
If $\Pi W\op \Pi=0$ then $z'=0$ and $y$ is constant at leading order.
\end{rem}

\subsection{Finite dimensional case study}
\label{Sec_5cas}

We now discuss the implications of Lemma~\ref{minispectre} for the different
values of $\mu/p$. We have already seen that a transition occurs in the 
definition of $B^\eps$ when $p=\mu$. Another transition happens when
$\nu=0$, \textit{i.e.} $\mu=2p$. The main results of this part are
summarized in Table~\ref{table1}.

\paragraph{Case when $0<\mu<p$.}

If $0<\mu<p$, we may apply Corollary~\ref{coro} with $\nu=\mu$ and the
limit operator $B^0$ is
\[
B^0(n,m) = C(n,m) \frac1{\gamma(n,m)} \un(\delta(n,m)\neq0).
\]
In this circumstance there is no need to separate the cases when 
$\delta(n,m)\neq0$, since
\[
(\eps^{-\mu} A + \eps^{-\nu}B^\eps)(n,m) 
= \eps^{-\mu} \left( \frac{C(n,m)}{\gamma(n,m)} + o(1) \right).
\]
Hence for large values of $p$ compared to $\mu$, the dynamics is the same for 
almost degenerate levels and for exactly degenerate levels. After a time-layer 
of size $O(\eps^\mu)$ the system is driven by Eq.~(\ref{ed_z}), namely
\[
\partial_t z= \Pi W \op \Pi z, \qquad z(0) = \Pi y(0),
\]
where $\Pi$ the orthogonal projection onto 
$\Ker A\op \cap \Ker B^0\op =\Ker \Psidom\op$.

\paragraph{Case when $\mu=p$.}

If $\mu=p$, we again apply Corollary~\ref{coro} with $\nu=\mu$, but the 
form of $B^0$ is slightly different, and now $\delta(n,m)$ plays a r\^ole, 
namely
\[
B^\eps(n,m)=B^0(n,m) = C(n,m) \frac{\gamma(n,m)}{\gamma(n,m)^2+\delta(n,m)^2}
\un(\delta(n,m)\neq0).
\]
The limit equation is as in the previous case.

\paragraph{Case when $p<\mu<2p$.}

If $p<\mu<2p$, we really have three different orders of magnitude in 
Eq.~(\ref{ed}). The form for $B^0$ is now
\[
B^0(n,m) = C(n,m) \frac{\gamma(n,m)}{\delta(n,m)^2} \un(\delta(n,m)\neq0),
\]
and the size of the time-layer is $O(\eps^{2p-\mu})$. The limit equation is 
still unchanged.

\paragraph{Case when $\mu=2p$.}

If $\mu=2p$, then $\nu=0$ and the contribution of $B^0$ competes with that of 
$W$. Therefore the projector mentioned in Corollary~\ref{coro} is here 
$\Pi_A$, the orthogonal projection on $\Ker A\op$, and the limit equation is
\[
\d_t z = \Pi_A (B^0+W)\op \Pi_A z, \qquad z(0) = \Pi_A y(0),
\]
where
\[
B^0(n,m) = C(n,m) \frac{\gamma(n,m)}{\delta(n,m)^2} \un(\delta(n,m)\neq0),
\]
after a time-layer of size $O(\eps^{\mu})$. In this case, and in the following 
one, the projection only depends on $A$ (and not on $B^0$). The 
constant $c$ is simply $c=\min_{\|x\|\leq 1, \Pi_A x = 0} - (A\op x,x)$.

\paragraph{Case when $\mu>2p$ or $B_\eps\equiv0$.}

If $\mu>2p$, then $\nu<0$ and $B^\eps\to0$. This case is therefore treated in 
the same way as when $B_\eps\equiv0$. In both cases, the limit equation is
\[
\d_t z = \Pi_A W\op \Pi_A z, \qquad z(0) = \Pi_A y(0),
\]
after a time-layer of size $O(\eps^{\mu})$. Hence for small values of $p$
compared to $\mu$, two almost degenerate levels $n$ and $m$ for which
$\delta(n,m)\neq0$ are already too far apart to resonate with the wave.

\subsection{Case when $\mu=0$.}
\label{Sec_5mu0}

In the case when $\mu=0$, Section~\ref{Sec_dom} does not yield an interesting
result since $\rhod-\rhodeux$ is of order $O(1)$. Therefore we have to use the
transition rates obtained in Section~\ref{Sec_ave}, namely
\[
\Psiave(n,m) = 2 |V(n,m)|^2 \sum_{\beta \in \Z^r} \frac{\gamma(n,m)}
{\gamma(n,m)^2 + |\omega(m,n)+ \beta\cdot\omega+\eps^p \delta(n,m)|^2}
|\phi_\beta|^2.
\]
For a finite or an infinite number of energy levels, we have at leading order 
$\rhod=\rhoapp$, where 
\[
\d_t \rhoapp = (\Psiavezero+W)\op \rhoapp
\]
and
\[
\Psiavezero(n,m) = 2 |V(n,m)|^2 \sum_{\beta \in \Z^r} \frac{\gamma(n,m)}
{\gamma(n,m)^2 + |\omega(m,n)+ \beta\cdot\omega|} |\phi_\beta|^2.
\]
The discussion of the long time behavior and equilibrium state for such an 
equation is the same as for Eq.~\eqref{rhod2} (see Appendix~\ref{Sec_relax}).

\subsection{Infinite dimensional case study}
\label{Sec_5infini}

In the infinite dimensional case the convergence of $B^\eps$ towards $B^0$ is 
not sufficient to conclude. However, we can use series expansions when they 
only include one non-positive order in $\eps$, \textit{i.e.} 
\[
\Psidom + W = \eps^{-\sigma} \tilde A + o(1).
\]
In this case, we consider at leading order the solution $\rhoapp$ of 
\begin{equation}
\label{rhodapp}
\d_t \rhoapp = \eps^{-\sigma} \Psiapp\op \rhoapp
\end{equation}
where $\Psiapp$ is homogeneous of order $O(1)$. In this case, time has to be 
changed into $t'=\eps^{-\nu}t$. With this time scale, the equation is not 
singular any more. There remains to list the cases when such an homogeneous
rate operator occurs, restricting the discussion to $\mu>0$. As for a finite 
number of levels, the value of the ratio $\mu/p$ is crucial. We once more have 
to consider the two cases $\mu\leq p$ and $\mu\geq p$. The main results of 
this section are summarized in Table~\ref{table2}.

\paragraph{Case when $\mu<p$.}
If $\mu<p$, series expansions yield that
\[
\Psidom(n,m) = C(n,m) \frac{\eps^{-\mu}}{\gamma(n,m)} + O(\eps^{2 p-3\mu}).
\]
If $2 p-3\mu>0$, i.e. $\mu < 2p/3$
\[
\Psidom(n,m) =: \eps^{-\mu} \Psiapp(n,m) + o(1).
\]
Since $\nu>0$, $\sigma=\mu$ and rates are homogeneous only if $W=0$. \\

\paragraph{Case when $\mu=p$.} 

In the case when $\mu=p$ no series expansion is needed and we have exactly
\[
\Psidom(n,m) 
= C(n,m) \frac{\eps^{-\mu}\gamma(n,m)}{\gamma(n,m)^2+\delta(n,m)^2} 
=: \eps^{-\mu} \Psiapp(n,m).
\]
Once more, we can only conclude if $W=0$.

\paragraph{Case when $\mu>p$.}

If $\mu<p$, series expansions lead to
\begin{eqnarray*}
\Psidom(n,m) & = & C(n,m) \frac{\eps^{-\mu}}{\gamma(n,m)} \un(\delta(n,m)=0) \\
& + & C(n,m) \frac{\eps^{-(2p-\mu)}\gamma(n,m)}{\delta(n,m)^2}
\un(\delta(n,m)\neq0) + O(\eps^{3\mu-4p}).
\end{eqnarray*}
Two cases lead to only one term $\eps^{-\sigma} \tilde A + o(1)$, either
$\delta(n,m)$ is always nonzero, $\mu > \frac43 p$ and $W=0$. In this case
$\sigma=2p-\mu$. Or $\mu > 2 p$ and $\sigma=\mu$. We once more need to assume 
that $W=0$. \\

In all the above cases, using the same type of estimates and integral
formulations as in the proofs in Section~\ref{Sec_ave} and \ref{Sec_dom}, we
define an approximate solution $\rhoapp$ to Eq.~\eqref{rhodapp} such that
\[
\|\rhod - \rhoapp\|_{L^\infty([0,T],l^2)} = o(1),
\]
and $\Psiapp$ has the properties described in Section~\ref{Sec_relax}.

\begin{table}[H]
\[
\begin{array}{c|c|c}
\mu/p & \sigma 
& \disp \frac{\Psiapp(n,m)}{2|V(n,m)|^2\sum_\beta|\phi_\beta|^2} \\
\hline \hline 0 \ (\mu=0,\ p>0) & \mu & \disp \frac1{\gamma(n,m)} 
\textrm{ (see caption)}
\\
\hline 0 < \mu/p < 2/3 & \mu & \disp \frac1{\gamma(n,m)} \\
\hline \mu/p = 1 
& \mu & \disp \frac{\gamma(n,m)}{\gamma(n,m)^2+\delta(n,m)^2} \\
\hline 4/3 < \mu/p < 2 
& 2p-\mu & \disp \frac{\gamma(n,m)}{\delta(n,m)^2} \\
\hline \mu/p = 2 & 0 & \disp \frac{\gamma(n,m)}{\delta(n,m)^2} \\
\hline 2 < \mu/p < \infty & 0 & 0 
\end{array}
\]
\caption{\label{table2}Cases when the transition rates are single powers of 
$\eps$, for an infinite number of levels. The sum over $\beta$'s includes only 
resonant contributions in general, except when $\mu=0$ where $\beta\in\Z^r$.}
\end{table}

\section{Restriction to a finite number of levels}
\label{Sec_fini}

It is important to show that we can restrict the study to a finite number of
levels since we have seen in Section~\ref{Sec_lim} that we are able to give
precise results on the time evolution of rate equations only in this case.
Another perspective is numerical simulations, which in any case can only treat
finite data.

Consider a solution $\rho$ to Bloch equations~\eqref{bloch} with initial datum
$\rho(0)$, and infinitely many quantum levels. In this section, we show how
$\rho$ may be approximated by $\rhoN$, solution to Eq.~\eqref{bloch} with
only a finite number $N$ of levels.

For this purpose, for all $N \in \N$, we define $\pi^N$, the projection of the
space $\C^{\N^2}$ of infinite matrices onto the space of $N \times N$
matrices, by
\[
(\pi^N u)(n,m) := u(n,m) \un[n,m < N].
\]
Then, the $N$-level \textit{truncated system}~(\ref{bloch}$^N$) is defined
from Eq.~\eqref{bloch} by
\begin{equation*}
\tag{\ref{bloch}$^N$}
\begin{split}
\eps^2 \d_t \rhoN(t,n,m)
& = -i \omega_\eps(n,m) \rhoN(t,n,m) + (\pi^N Q_\eps)(\rhoN)(n,m)
\\
& + i \eps \sum_{k<N}
\left[
{\cal V}\left(\frac{t}{\eps^2},n,k\right) \rhoN(t,k,m)
-
{\cal V}\left(\frac{t}{\eps^2},k,m\right) \rhoN(t,n,k)
\right].
\end{split}
\end{equation*}
The initial datum $\rhoN(0)$ is also naturally defined as
\[
\rhoN(0) := \pi^N \rho(0).
\]

The analysis of the previous sections shows that, in the limit $\eps\to0$,
$\rho$ is approximated by the diagonal solution $\rhodeux$ to rate
equations~\eqref{rhod2}. Theorem~\ref{Th_relax} also gives an approximation of
$\rhoN$ by $\rhoNdeux$, which turns out to be solution to the truncated system
obtained from Eq.~\eqref{rhod2}:
\begin{equation*}
\tag{\ref{rhod2}$^N$}
\d_t \rhoNdeux = (\pi^N \Psidom + \pi^N W)\op \rhoNdeux , \quad
\rhoNdeux(0) = \pi^N \rho(0).
\end{equation*}

To apply Theorem~\ref{Th_relax}, we need Hypotheses~\ref{dioph} to \ref{far}.
Under the additional condition that longitudinal relaxation coefficients are
decaying enough at infinity (Hypothesis~\ref{Wdecay}), we show that this
truncation procedure is compatible with the evolution according to Bloch
equations and rate equations:

\begin{lemma}
Under Hypotheses~\ref{dioph} to \ref{Wdecay}, for all $\nu, T, \eps > 0$,
there exists an integer $N$ such that,
\[
\textrm{if } \quad \| (1-\pi^N)\rho(0) \|_{l^2} \leq \nu, \textrm{ then }
\| \rho - \rhoN \|_{L^\infty([0,T],l^2)}
\leq 2 \nu + C (\eps^\mu+\eps^{1-2\mu}),
\]
where $C=C(T)$ is the constant from Theorem~\ref{Th_relax}, and $N$ has the
form
\[
N=N_0(\nu,T) \eps^{-\frac{\mu}{(1+\eta)N_\eta}},
\]
with $\eta$ and $N_\eta$ given by Hypotheses~\ref{dioph} and \ref{extrasmooth}.

If, in addition, no resonance occurs between the wave and high energy levels
(\textit{i.e.} there exists $M \in \N$ such that, when $\min(n,k)>M$, the set
$\{ \beta \in \Z^r \suchthat \omega(n,k) + \beta\cdot\omega = 0 \}$ is empty), 
then $N$ has the form $N=N_0(\nu,T)$, uniformly with respect to $\eps$.
\end{lemma}

\begin{proof}
Theorem~\ref{Th_relax} applies for both the infinite and the finite number of 
levels problems, therefore
\[
\|\rho-\rhodeux \|_{L^\infty([0,T],l^2)} \leq C (\eps^\mu + \eps^{1-2\mu})
\]
and
\[
\|\rhoN-\rhoNdeux \|_{L^\infty([0,T],l^2)} \leq C (\eps^\mu + \eps^{1-2\mu}),
\]
and we only need to estimate the difference $\Delta := \rhodeux-\rhoNdeux$, 
which is solution to
\[
\d_t \Delta = \left(\Psidom + W\right)\op \Delta 
+ \left(\Psidom - \PsidomN\right)\op \rhoNdeux
+ \left(W - \pi^N W\right)\op \rhoNdeux.
\]
Thanks to the non-positiveness property of the operators associated with
$\Psidom$ and $W$, an integral formulation leads to
\begin{eqnarray*}
\|\Delta\|_{L^\infty([0,T],l^2)}
& \leq & C \Big(
\|(1-\pi^N)\rho(0)\|_{l^2}
+ \|\left(\Psidom - \pi^N\Psidom\right)\op \rhoNdeux\|_{L^\infty([0,T],l^2)} \\
&& \qquad + \|\left(W - \pi^N W\right)\op \rhoNdeux\|_{L^\infty([0,T],l^2)} 
\Big) \\
& \leq &
C \Big( \|(1-\pi^N)\rho(0)\|_{l^2}
+ \|\Psidom - \pi^N\Psidom\|_{l^\infty_k l^1_n \cap l^\infty_n l^1_k} \\
&& \qquad + \|W - \pi^N W\|_{l^\infty_k l^1_n \cap l^\infty_n l^1_k}
\Big).
\end{eqnarray*}

The first term goes to zero as $N$ goes to infinity simply because the initial
datum $\rho(0)$ is in $l^2$. The third one reads
\begin{multline*}
\sup_{n>N} \sum_k |W(n,k)| + \sup_n \sum_{k>N} |W(n,k)|
+ \sup_{k>N} \sum_n |W(n,k)| + \sup_k \sum_{n>N} |W(n,k)| \\
\begin{aligned}
& \leq C \left( \sup_{n>N} \sum_k |W(n,k)| + \sup_n \sum_{k>N} |W(n,k)| \right)
 \\
& \leq C N^{-K} \left( \sup_{n>N} \sum_k (1+n)^K |W(n,k)|
+ \sup_n \sum_{k>N} (1+k)^K |W(n,k)| \right) \\
& \leq C N^{-K},
\end{aligned}
\end{multline*}
thanks to Hypothesis~\ref{Wdecay}. Thus, this term is also $o(1)$ uniformly
with respect to $\eps$ as $N$ goes to infinity.

Finally, taking into account the fact that $\Psidom$ is symmetric, the second
term is
\begin{multline*}
2 \left( \sup_{n>N} \sum_k |\Psidom(n,k)| 
+ \sup_n \sum_{k>N} |\Psidom(n,k)| \right) \\
\begin{aligned}
&\leq C \eps^{-\mu} \sup_{n>N} \sum_k
\sum_{\beta \suchthat \omega(n,k)+\beta\cdot\omega = 0} 
|\phi_\beta|^2 |V(n,k)|^2 \\
&\leq C \eps^{-\mu} N^{-(1+\eta)N_\eta} \sup_{n>N} 
\sum_k (1+n)^{(1+\eta)N_\eta} 
\sum_{\beta \suchthat \omega(n,k)+\beta\cdot\omega = 0} 
|\phi_\beta|^2 |V(n,k)|^2 \\
&\leq C \eps^{-\mu} N^{-(1+\eta)N_\eta},
\end{aligned}
\end{multline*}
which vanishes in fact when no resonance occurs between the wave and high
energy levels. Else, we obtain a $o(1)$ as $\eps$ goes to zero under the
condition
\[
N \gg \eps^{-\mu/N_\eta(1+\eta)}.
\]
\end{proof}

\section{Appendix}
\label{Sec_app}

In this section, we give the lemmas concerning the modified relaxation 
operator from Eq.~\eqref{rate}, in order to describe the dynamics of the 
solution. We also give a proof of the genericity of the small divisor 
estimates of Hypothesis~\ref{dioph}. 

\subsection{Relaxation operators}
\label{Sec_relax}

We first prove some non-positiveness properties of $(\Wmod)\op$.

\subsubsection{Continuity and non-positiveness}

\begin{lemma} 
\label{nonpositif}
Let $A(n,m) \in l^\infty_n l^1_m \cap l^\infty_m l^1_n$. 

(i) Its associated operator $A\op$ is bounded on the spaces $l^q$, 
$1 \leq q \leq \infty$, and 
\[
\|A\op u\|_{l^q} \leq \|A(n,m)\|_{l^\infty_n l^1_m \cap l^\infty_m l^1_n} 
\|u\|_{l^q}.
\] 

(ii) If in addition $A(n,m) \geq 0$, then for all positive integer $N$, the 
spectrum of the restriction of $A\op$ to $\R^N$ is contained in
$\{ \re \lambda < 0 \} \cup \{ 0 \}$. 

(iii) If $A(n,m) \geq 0$ is symmetric, $A\op$ is non-positive on $l^2$, and 
the exponential $\exp(tA\op)$ is well defined as an operator on $l^2$ when 
$t\geq0$. Its norm is 1.
\end{lemma}

\begin{proof}

\textit{(i)} For $q=\infty$ or $q=1$, the result is immediate. The remaining 
cases are obtained by interpolation.

\textit{(ii)} The localization of the eigenvalues of $A\op$ is obtained 
\textit{via} the Hadamard-Gerschgorin method applied to 
$M := \transp \overline{A\op} = \transp A\op$, whose eigenvalues are the
conjugates of those of $A\op$: if $\lambda$ is an eigenvalue of $M$, there 
exists an index $n$ such that
\[
|\lambda-M(n,n)| \leq \sum_{m \neq n} |M(n,m)|.
\]
Remarking that for $m \neq n$, $M(n,m)=A\op(m,n)\geq0$, and 
$M(n,n)=-\sum_{m \neq n}M(n,m)$, the conclusion is straightforward.

\textit{(iii)} In the symmetric case, compute for all $u \in l^2$,
\begin{eqnarray*}
(A\op u,u) 
& = & \sum_n \sum_{m \neq n} A\op(n,m) u(m) u(n) +  \sum_n A\op(n,n) u(n)^2 \\
& \leq & \frac12 \sum_n \sum_{m \neq n} A\op(n,m) (|u(m)|^2 + |u(n)|^2) 
- \sum_n \sum_{m \neq n} A\op(m,n) u(n)^2 = 0.
\end{eqnarray*}
The norm of $\exp(tA\op)$ is 1 on $l^2$ for $t>0$ because $A\op$ has a 
non-trivial kernel (see below).
\end{proof}

In the same way as for item (i) of Lemma~\ref{nonpositif}, we have the 
estimate usually called Schur's lemma.

\begin{lemma}
\label{schur}
Let $A\in l_n^\infty l_k^1\cap l_k^\infty l_n^1$ and $1 \leq q \leq \infty$. 
Then, if $u\in l^q_{n,m}(\N\times\N)$, 
\[
\left\| \sum_k A(n,k) u(k,m) - A(k,m) u(n,k) \right\|_{l^q_{n,m}}
\leq
\|A(n,m)\|_{l^\infty_n l^1_m \cap l^\infty_m l^1_n} \|u(n,m)\|_{l^q_{n,m}}.
\]
\end{lemma}

\subsubsection{Asymptotic states of the rate equation}
\label{Sec_astate}

We define the asymptotic state $\rhoe$ associated with Eq.~\eqref{rate} and 
with the initial datum $\rho(0)$, the limit (in $l^1$), if it exists, of the 
solution $\rhod(t)$ to Eq.~\eqref{rate} with the initial value $\rho(0)$ as 
$t$ goes to infinity. Such an asymptotic state is necessarily an equilibrium 
state (\textit{i.e.} it belongs to the kernel of $(\Wmod)\op$). As an example, 
the usual thermodynamic equilibrium $\rhoe^{\rm therm}$ (corresponding to 
$\Wmod = W$ satisfying Eq.~\eqref{thermo}) for the $N$-level Bloch model is 
given by 
\[
\rhoe^{\rm therm}(n) =
\frac{\exp \left( -\frac{\omega_n}{T} \right) }
{\sum_{k=1}^N \exp \left( -\frac{\omega_k}{T} \right) }.
\]

We study the kernel of operators $A\op$ modelled on $(\Wmod)\op$: the
--finite or infinite-- matrix $A$, written in the eigenstates basis
$e=(e_1,e_2,\dots)$, has the property (${\cal P}$).
\begin{equation}
A(m,n)=0 \Leftrightarrow A(n,m)=0. \tag{$\cal P$}
\end{equation}
In particular, thanks to this property, a vanishing column in $A\op$ 
corresponds to a vanishing line, and conversely.

The kernel of $A\op$ is linked to the spaces generated by elements of the 
basis $e$,
\[
E_L := {\rm Span} \{ e_m \suchthat m \in L \}, 
\]
when $L$ is a finite subset of indices.
 
\begin{proposition}
\label{Pr_noyau1d}
Let $A(m,n) \in l^\infty_m l^1_n \cap l^\infty_n l^1_m$ satisfy property 
(${\cal P}$), and $A(m,n) \geq 0$ when $m \neq n$. In addition, suppose that 
there exists a decomposition of $l^1$ into $A\op$-stable subspaces of the form 
$E_L$ (each with finite dimension). \\
Then, the restriction of $A\op$ to any such non-zero subspace with minimal 
size has a one-dimensional kernel.
\end{proposition}

Such a minimal subset can be one dimensional, \textit{i.e.} generated by a 
single $e_n$. This corresponds exactly to the case when the $n$-th line (or 
column) of $A\op$ vanishes. When the cardinality of $L$ is greater or equal to 
two, $E_L$ is characterized by:
\[
\begin{aligned}[t]
m \leq n \in L \quad \Leftrightarrow \quad \exists m=:m_1 & \leq \dots \leq
m_s:=n, \textrm{ such that} \\
& m_j \in L \textrm{ for all } j=1,\dots,s \\
\mbox{and } & A\op(m_j,m_{j+1}) \neq 0 \textrm{ for all } j=1,\dots,s-1.
\end{aligned}
\]

\begin{proof}
When the dimension of $E_L$ is one, the result is trivial. Now, suppose that 
the cardinality of $L$ is at least $2$. 

Denote by $A\op^L$ the restriction of $A\op$ to $E_L$. Since the dimension of 
$E_L$ is finite, the dimension of the kernel of $A\op^L$ is the same as the 
dimension of the kernel of the transposed matrix, $\transp A\op^L$. Let $u$ 
belong to this kernel. Written relatively to the basis 
$\{ e_n \suchthat n  \in L \}$, the relation $\transp A\op^L u=0$ reads
\[
\left( \sum_{k \in L \setminus \{n\}} A(n,k) \right) u(n) = 
\sum_{k \in L \setminus \{n\}} A(n,k) u(k)
\]
for all $n \in L$. Each coefficient $A(n,k)$ is non-negative, and thanks to 
the property (${\mathcal P}$), the sum $\sum_{k \in L \setminus \{n\}} A(n,k)$ 
is positive. Therefore, $u(n)$ belongs to the convex hull of the other 
coordinates. Since this is valid for all $n \in L$, all the coordinates must 
be equal, and the kernel of $\transp A\op^L$ is generated by 
$\sum_{n \in L} e_n$.
\end{proof}

\begin{Corollary}
Under the assumptions of Proposition~\ref{Pr_noyau1d}, for each initial datum 
$\rho(0) \in l^1$ satisfying Eq.~\eqref{datum}, there is a unique 
$\rhoe \in l^1 \cap \Ker A\op$ such that for all minimal $A\op$-stable 
subspace $E_L$,  
\[
\sum_{n \in L} \rhoe(n) = \sum_{n \in L} \rho(0,n).
\]
\end{Corollary}

Since the decomposition of $l^1$ into $A\op$-stable finite dimensional  
subspaces corresponds to the splitting of Eq.~\eqref{rate} into a countable 
set of decoupled finite dimensional systems which have bounded solutions 
(thanks to the results of Lemma~\ref{nonpositif} on the spectrum of $A\op$), a 
diagonal argument allows to extract a converging subsequence of $\rhod(t)$. 
Finally, the normalization of the trace on each minimal subspace ensures the
uniqueness of the limit. This shows that the whole sequence 
$(\rhod(t))_{t \geq 0}$ converges to $\rhoe$, which is actually an asymptotic 
state.

In the purely infinite dimensional case, \textit{i.e.} when no
decomposition into finite dimensional $A\op$-stable subspaces exists,
the results above may break down.
\begin{proposition}
\label{Pr_noyau0d}
Consider $A(m,n) \in l^\infty_m l^1_n \cap l^\infty_n l^1_m$ satisfying
property (${\cal P}$), $A(m,n) \geq 0$ when $m \neq n$, and either symmetric, 
or in Pauli form (relation~\eqref{thermo}). Suppose that there exists a 
minimal $A\op$-stable subspace $E$ of $l^1$ generated by an infinite number of 
eigenstates $e_n$. 

Then, the kernel of ${A\op}_{|_{E}}$ (the restriction to $E$ of the operator 
on $l^1$) is $\{ 0 \}$.
\end{proposition}

\begin{proof}
Denoting $A\op^E$ the restriction of $A\op$ to $E$, we have a one-to-one 
relation between elements of the kernel of $A\op^E$ and elements of the kernel 
of $\transp A\op^E$:
\begin{align*}
u \in \Ker A\op^E 
& \Leftrightarrow \transp u \in \Ker \transp A\op^E 
&& \textrm{ in the symmetric case}, \\
& \Leftrightarrow \left( u(n) \exp\left( \frac{\omega(n)}{T} \right) 
\suchthat n \in L \right) \in \Ker \transp A\op^E 
&& \textrm{ in the ``Pauli'' case}.
\end{align*}
Since we have the bounds 
\[
0 < \exp\left( \frac{\omega(1)}{T} \right) 
\leq \exp\left( \frac{\omega(n)}{T} \right) \leq 
\exp\left( \frac{\omega_{\rm ionisation}}{T} \right)
\]
for all $n$, this correspondence preserves the summability property.

Finally, the proof of Proposition~\ref{Pr_noyau1d} shows that the kernel of 
$\transp A\op^E$ (in $l^1$) is $\{ 0 \}$, and this gives the result.
\end{proof}

\begin{Corollary}
\label{noeq}
Under the assumptions of Proposition~\ref{Pr_noyau0d}, for any initial datum 
$\rho(0)$ with non-vanishing component in $E$, there is no equilibrium state
$\rhoe \in l^1 \cap \Ker A\op$ with the same trace as $\rho(0)$ in $E$.
\end{Corollary}

\begin{rem}
In the symmetric case (when $\Psidom=0$, $\Wmod=W$), according to
Lemma~\ref{nonpositif}, the $l^2$-norm of the solution $\rhod$ to 
Eq.~\eqref{rate} is decreasing in time; thus, it tends to a certain
value $r \geq 0$. This means that $\rhod$ approaches a limit cycle in $l^2$
belonging the intersection of the sphere $\|\rhod\|_{l^2} = r$ and the
hyperplane where the $l^1$-norm is one (assuming for simplicity that
there is no strict $A\op$-stable subspace of $l^1$). In this case,
only weak convergence (to zero) can occur. 
\end{rem}

\subsection{Diophantine estimates}
\label{Sec_dioph}

We show the genericity of Hypothesis~\ref{dioph}.

\begin{lemma}
For all $\eta>0$ and all real sequence $\omega(n,m)$, there exists a constant 
$C_\eta>0$, such that for almost all value of the frequency vector
$\omega=(\omega_1,\dots,\omega_r)$,
\begin{eqnarray*}
&&
\non
\forall \alpha=(\alpha_1,\ldots,\alpha_r) \in \Z^r \setminus \{ 0 \}, \quad
\forall (n,k) \in \N^2 \text{ such that }
\alpha \cdot \omega + \omega(n,k) \neq 0,
\\
&&
\qquad
|\alpha \cdot \omega+\omega(n,k)| \geq
\frac{C_\eta}{(1+|\alpha|)^{r-1+\eta} (1+n)^{1+\eta} (1+k)^{1+\eta}}.
\end{eqnarray*}
\end{lemma}

\begin{proof}
We follow the standard approach (see \textit{e.g.}
\cite{AG}). Restricting $\omega$ to a ball $B$ in $\R^r$, we show that
the measure of the set of ``bad frequencies'' violating the inequality
for all constant $C$ is zero. 

For $\eta,c>0$, $\alpha \in \Z^r \setminus \{ 0 \}$ and $(n,k) \in \N^2$ fixed,
set
\[
B_{\alpha,n,k}^{\eta,c} 
:= \left\{ \omega \in B \suchthat |\alpha \cdot \omega+\omega(n,k)| \leq
\frac{c}{(1+|\alpha|)^{r-1+\eta} (1+n)^{1+\eta} (1+k)^{1+\eta}} \right\}.
\]
This limitates $\omega$ in the direction of $\alpha$. Introducing a
constant $K$ which depends on the size of $B$ only, we obtain
\[
\textrm{meas} \left( B_{\alpha,n,k}^{\eta,c} \right) 
\leq \frac{Kc}{(1+|\alpha|)^{r-1+\eta} (1+n)^{1+\eta} (1+k)^{1+\eta}}.
\]
Now, with $\eta,c>0$ fixed, the measure of the set of frequencies for which 
the inequality is false at least for some $(\alpha,n,k)$ is less than the sum 
(over $\alpha$, $n$ and $k$) of the ones above, and thus is $O(c)$.
\end{proof}

\medskip

{\bf Acknowledgments:}
This work has been partially supported by the GDR ``Amplitude Equations and 
Qualitative Properties'' (GDR CNRS 2103 : EAPQ) and the European Program 
'Improving the Human Potential' in the framework of the 'HYKE' network
HPRN-CT-2002-00282. The authors thank Laurent Bonavero for fruitful 
discussions.


\newcommand{\etalchar}[1]{$^{#1}$}


\begin{thebibliography}{CTDRG88}

\bibitem[AG91]{AG}
S.~Alinhac and P.~G\'erard.
\newblock {\em Op\'erateurs pseudo-diff\'erentiels et th\'eor\`eme de
  Nash-Moser}.
\newblock Inter-Editions, 1991.

\bibitem[Arn89]{Ar}
V.I. Arnol'd.
\newblock {\em Mathematical Methods of Classical Mechanics}.
\newblock Number~60 in Graduate Texts in Mathematics. Springer-Verlag, 1989.

\bibitem[BBR01]{BBR}
B.~Bid\'egaray, A.~Bourgeade, and D.~Reignier.
\newblock Introducing physical relaxation terms in {B}loch equations.
\newblock {\em J. Comput. Phys.}, 170(2):603--613, 2001.

\bibitem[BCEP03]{BCEP}
D.~Benedetto, F.~Castella, R.~Esposito, and M.~Pulvirenti.
\newblock Some considerations on the derivation of the nonlinear quantum
  {B}oltzmann equation.
\newblock {\em J. Stat. Phys.}, 2003.
\newblock To appear.

\bibitem[BFCD03]{BCD}
B.~Bid\'egaray-Fesquet, F.~Castella, and P.~Degond.
\newblock From {B}loch model to the rate equations.
\newblock {\em Discrete Contin. Dynam. Syst.}, 2003.
\newblock To appear.

\bibitem[{Bid}03]{Bi}
B.~{Bid\'egaray-Fesquet}.
\newblock {\em De Maxwell-Bloch \`a Schr\"odinger non lin\'eaire~: une
  hi\'erarchie de mod\`eles en optique quantique}.
\newblock 2003.
\newblock In preparation.

\bibitem[Boh79]{Boh}
A.~Bohm.
\newblock {\em Quantum Mechanics}.
\newblock Texts and monographs in Physics. Springer-Verlag, 1979.

\bibitem[Boy92]{Boy}
R.W. Boyd.
\newblock {\em Nonlinear Optics}.
\newblock Academic Press, 1992.

\bibitem[Cas99]{Ca1}
F.~Castella.
\newblock On the derivation of a quantum {B}oltzmann equation from the periodic
  von {N}eumann equation.
\newblock {\em M2AN}, 33(2):329--349, 1999.

\bibitem[Cas01]{Ca3}
F.~Castella.
\newblock From the von {N}eumann equation to the quantum {B}oltzmann equation
  in a deterministic framework.
\newblock {\em J. Stat. Phys.}, 104(1/2):387--447, 2001.

\bibitem[Cas02]{Ca2}
F.~Castella.
\newblock From the von {N}eumann equation to the quantum {B}oltzmann equation
  {II}: identifying the {B}orn series.
\newblock {\em J. Stat. Phys.}, 106(5/6):1197--1220, 2002.

\bibitem[CCC{\etalchar{+}}03]{CCCFLLT}
E.~Canc\`es, F.~Castella, P.~Chartier, E.~Faou, C.~{Le Bris}, F.~Legoll, and
  G.~Turinici.
\newblock Long-time averaging using symplectic solvers with applications to
  molecular dynamics.
\newblock 2003.
\newblock In preparation.

\bibitem[CP02]{CP1}
F.~Castella and A.~Plagne.
\newblock A distribution result for slices of sums of squares.
\newblock {\em Math. Proc. Cambridge Philos. Soc.}, 132(1):1--22, 2002.

\bibitem[CP03]{CP2}
F.~Castella and A.~Plagne.
\newblock Non-derivation of the quantum {B}oltzmann equation from the periodic
  {S}chr\"odinger equation.
\newblock {\em Indiana Univ. Math. J.}, 51(4):963--1016, 2003.

\bibitem[CTDRG88]{CTDRG}
C.~Cohen-Tannoudji, J.~Dupont-Roc, and G.~Grynberg.
\newblock {\em Processus d'interaction entre photons et atomes}.
\newblock Savoirs actuels. Intereditions/Editions du CNRS, 1988.

\bibitem[EY00]{EY}
L.~Erd\"os and H.T. Yau.
\newblock Linear {B}oltzmann equation as the weak coupling limit of a random
  {S}chr\"odinger equation.
\newblock {\em Comm. Pure Appl. Math.}, 53(6):667--735, 2000.

\bibitem[KL57]{KL1}
W.~Kohn and J.M. Luttinger.
\newblock Quantum theory of electrical transport phenomena.
\newblock {\em Phys. Rev.}, 108(3):590--611, 1957.

\bibitem[KPR96]{KPR}
J.B. Keller, G.~Papanicolaou, and L.~Ryzhik.
\newblock Transport equations for elastic and other waves in random media.
\newblock {\em Wave Motion}, 24(4):327--370, 1996.

\bibitem[Kre83]{Kr}
H.J. Kreuzer.
\newblock {\em Nonequilibrium thermodynamics and its statistical foundations}.
\newblock Monographs on Physics and Chemistry of Materials. Oxford Science
  Publications, 1983.

\bibitem[Lin76]{Li}
G.~Lindblad.
\newblock On the generators of quantum dynamical semigroups.
\newblock {\em Comm. Math. Phys.}, 48:119--130, 1976.

\bibitem[LK58]{KL2}
J.M. Luttinger and W.~Kohn.
\newblock Quantum theory of electrical transport phenomena. {II}.
\newblock {\em Phys. Rev.}, 109(6):1892--1909, 1958.

\bibitem[LM88]{LM}
P.~Lochak and C.~Meunier.
\newblock {\em Multiphase averaging for classical systems. With applications to
  adiabatic theorems}.
\newblock Number~72 in Applied Mathematical Sciences. Springer-Verlag, 1988.

\bibitem[Lou91]{Lo}
R.~Loudon.
\newblock {\em The quantum theory of light}.
\newblock Clarendon Press, Oxford, 1991.

\bibitem[Nie96]{Ni}
F.~Nier.
\newblock A semi-classical picture of quantum scattering.
\newblock {\em Ann. Sci. Ec. Norm. Sup., 4. S\'er.}, 29(2):149--183, 1996.

\bibitem[NM92]{NM}
A.C. Newell and J.V. Moloney.
\newblock {\em Nonlinear Optics}.
\newblock Advanced Topics in the Interdisciplinary Mathematical Sciences.
  Addison-Wesley Publishing Company, 1992.

\bibitem[Spo77]{Sp1}
H.~Spohn.
\newblock Derivation of the transport equation for electrons moving through
  random impurities.
\newblock {\em J. Stat. Phys.}, 17(6):385--412, 1977.

\bibitem[Spo80]{Sp2}
H.~Spohn.
\newblock Kinetic equations from {H}amiltonian dynamics: {M}arkovian limits.
\newblock {\em Rev. Mod. Phys.}, 52(3):569--615, 1980.

\bibitem[Spo91]{Sp3}
H.~Spohn.
\newblock {\em Large Scale Dynamics of interacting particles}.
\newblock Texts and Monographs in Physics. Springer, Berlin, 1991.

\bibitem[SSL77]{SSL}
M.~Sargent, M.O. Scully, and W.E. Lamb.
\newblock {\em Laser Physics}.
\newblock Addison-Wesley, 1977.

\bibitem[SV85]{SV}
J.A. Sanders and F.~Verhulst.
\newblock {\em Averaging methods in nonlinear dynamical systems}, volume~59 of
  {\em Applied Mathematical Sciences}.
\newblock Springer-Verlag, 1985.

\bibitem[vH55]{VH1}
L.~van Hove.
\newblock Quantum-mechanical perturbations giving rise to a statistical
  transport equation.
\newblock {\em Phyisca}, 21:517--530, 1955.

\bibitem[vH57]{VH2}
L.~van Hove.
\newblock The approach to equilibrium in quantum statistics. {A} perturbation
  treatment to general order.
\newblock {\em Phyisca}, 23:441--480, 1957.

\bibitem[Zwa66]{Zw}
R.~Zwanzig.
\newblock {\em Quantum Statistical Mechanics}.
\newblock Gordon and Breach, New-York, 1966.

\end{thebibliography}
\end{document}